\documentclass[oneside]{amsart}
\usepackage{amssymb,latexsym,amsmath,amsthm}
\usepackage{fancyhdr}

\newtheorem{thm}{Theorem}
\newtheorem*{thm*}{Theorem}
\newtheorem{prop}{Proposition}
\newtheorem*{prop*}{Proposition}
\newtheorem{lemma}[prop]{Lemma}
\newtheorem*{lemma*}{Lemma}
\newtheorem{cor}{Corollary}[thm]
\newtheorem*{cor*}{Corollary}
\newtheorem*{conjecture*}{Conjecture}

\theoremstyle{definition}

\newcommand{\con}{\equiv}

\newcommand{\modd}[1]{\; ( \text{mod} \; #1)}
\newcommand{\bstack}[2]{\substack{#1 \\ #2}}

\newcommand{\maps}{\rightarrow}

\newcommand{\Union}{\bigcup}

\newcommand{\al}{\alpha}
\newcommand{\be}{\beta}
\newcommand{\gam}{\gamma}
\newcommand{\ga}{\gamma}
\newcommand{\del}{\delta}
\newcommand{\ep}{\epsilon}

\newcommand{\lam}{\lambda}

\newcommand{\Ga}{\Gamma}

\newcommand{\uchi}{\underline{\chi}}
\newcommand{\tchi}{\tilde{\chi}}

\newcommand{\A}{\mathcal{A}}
\newcommand{\Acal}{\mathcal{A}}
\newcommand{\Bcal}{\mathcal{B}}

\newcommand{\Ecal}{\mathcal{E}}

\newcommand{\Ical}{\mathcal{I}}

\newcommand{\Vcal}{\mathcal{V}}
\newcommand{\V}{\mathcal{V}}

\newcommand{\C}{\mathbb{C}}

\newcommand{\R}{\mathbb{R}}

\newcommand{\Z}{\mathbb{Z}}

\newcommand{\M}{\mathfrak{M}}
\newcommand{\m}{\mathfrak{m}}

\newcommand{\beq}{\begin{equation}}
\newcommand{\eeq}{\end{equation}}

\begin{document}

\title{Discrete fractional Radon transforms and quadratic forms}
\author{Lillian B. Pierce}
\address{School of Mathematics, Institute for Advanced Study, Princeton New Jersey 08540}
\curraddr{Mathematical Institute, University of Oxford, 24-29 St. Giles', Oxford OX1 3LB, United Kingdom }
\email{lillian.pierce@maths.ox.ac.uk}
\thanks{During this work, the author was supported by the Simonyi Fund at the Institute for Advanced Study and the National Science Foundation, including DMS-0902658 and DMS-0635607.}
\date{September 2012. Note: This is a corrected version of the original paper, which appeared in the {\emph{Duke Math. Journal}} {\bf{161}} no. 1 (2012) 69-106. The statements of Propositions 3, 6, 7, and Theorem 1 have been corrected, and Corollary 1.1 has been added.}
\maketitle

\begin{abstract}
We consider discrete analogues of fractional Radon transforms involving integration over paraboloids defined by positive definite quadratic forms. We prove that such discrete operators extend to bounded operators from $\ell^p$ to $\ell^q$ for a certain family of kernels. The method involves an intricate spectral decomposition according to major and minor arcs, motivated by ideas from the circle method of Hardy and Littlewood. Techniques from harmonic analysis, in particular Fourier transform methods and oscillatory integrals, as well as the number theoretic structure of quadratic forms, exponential sums, and theta functions, play key roles in the proof.

\end{abstract}

\section{Introduction}
The focus of this paper is a family of discrete Radon transforms along paraboloids, with singularities of subcritical degree. These belong to the world of fractional integral operators, which in the classical  setting of $\R^k$ are defined most simply by
\beq\label{simple_frxnl_cts}
\mathcal{I}_\lam f(x) = \int_{\R^k} \frac{f(x-y)}{|y|^{k\lam}} dy,
\eeq
with $0< \lam < 1$. The Hardy-Littlewood-Sobolev theorem states that when $1< p < q < \infty$ and $1/q = 1/p -(1-\lam)$, $\Ical_\lam$ extends to a bounded operator from $L^p(\R^k)$ to $L^q(\R^k)$.
The discrete analogue of this operator may be defined for  a (compactly supported) function $f: \Z^k \maps \C$ by
\beq\label{simple_frxnl}
 I_\lam f(n) = \sum_{\bstack{m \in \Z^k}{m \neq 0}} \frac{f(n-m)}{|m|^{k\lam}},
 \eeq
where $n \in \Z^k$ and $0< \lam < 1$. An elementary comparison argument shows that as a  consequence of the boundedness of the original operator $\Ical_\lam$, the discrete operator $I_\lam$ extends to a bounded operator from $\ell^p(\Z^k)$ to $\ell^q(\Z^k)$  for all $1/q \leq 1/p - (1-\lam)$ with $1 < p < q < \infty$; moreover, this result is sharp (see Proposition \emph{a} of \cite{SW2}).

Simple comparison arguments are no longer as useful for variants of $I_\lam$ of the form 
\beq\label{IP_lam}
 I^P_\lam f(n) = \sum_{\bstack{m \in \Z^{k_1}}{m \neq 0}}\frac{f(n-P(m))}{|m|^{k_1\lam}},
 \eeq
where $n \in \Z^{k_2}$, $0< \lam < 1$, and $P$ is a polynomial map from $\Z^{k_1}$ to $\Z^{k_2}$; this is now a discrete analogue of a fractional Radon transform. If $P$ is a one-to-one mapping, comparing $I^P_\lam$ to its continuous analogue again shows that $I^P_{\lam}$ maps $\ell^p(\Z^{k_2})$ to $\ell^q(\Z^{k_2})$ boundedly for $1/q \leq 1/p - (1-\lam)$ with $1 < p < q < \infty$ (see Proposition \emph{b} of \cite{SW2}). 

That this is not in general sharp may already be seen in the one-dimensional case $k_1=k_2 = 1$ with $P(m) = m^s$, $s \geq 2$ an integer. Then the operator
\beq\label{dfn_IPlam}
 I^s_\lam f(n) = \sum_{m=1}^\infty \frac{f(n-m^s)}{m^{\lam}}
 \eeq
is expected to be bounded from $\ell^p(\Z)$ to $\ell^q(\Z)$ whenever $p,q$ satisfy the two conditions
\newline
\indent (i) $1/q \leq 1/p  - (1-\lam)/s$,
\newline
\indent (ii) $1/q < \lam, 1/p > 1-\lam$.
\newline
This is only known in full for $s=2$, due to work of Stein and Wainger \cite{SW2}, \cite{SW3}, Oberlin \cite{Obe}, and Ionescu and Wainger \cite{IW}. For $s \geq 3$ only significantly weaker results are known for $I^s_\lam$, and its behavior is closely linked to major unsolved problems in number theory (precisely, Hypothesis $K^*$ of Hardy and Littlewood and certain aspects of Waring's problem; see \cite{SW2} and \cite{Pie_Waring_2010}).
It is the quadratic nature of the argument in $I^2_\lam$ (and hence of its Fourier multiplier) that makes this operator approachable while higher degree operators are not. 

\subsection{The main theorem}
In this paper we prove results for a broad class of discrete fractional Radon transforms of the form (\ref{IP_lam}) that retain the key quadratic nature of $I^2_\lam$ while introducing new number theoretic features and expanding the underlying space to higher dimensions.

 Let $Q_1, Q_2$ be positive definite quadratic forms in $k$ variables with integer coefficients, and define  for any (compactly supported) function $f:\Z^{k+1} \maps \C$ the operator
\beq\label{dis_JQ_dfn} 
 J_{Q_1,Q_2,\lam} f(n,t) = \sum_{\bstack{m \in \Z^k}{m \neq 0}} \frac{f(n-m,t-Q_1(m))}{Q_2(m)^{k\lam/2}},
 \eeq
 where $n \in\Z^k, t \in \Z$, and $0< \lam < 1$.

We prove the following theorem for the operator $J_{Q_1,Q_2,\lam}$:
\begin{thm}\label{dis_JQ_thm}
Let $\lam_k =0$ for $k=1,2$ and $\lam_k = k/(2k+2)$ for $k \geq 3$. For any $k \geq 1$ and $\lam_k< \lam < 1$, the operator $J_{Q_1,Q_2,\lam}$ extends to a bounded operator from $\ell^p(\Z^{k+1})$ to $\ell^q(\Z^{k+1})$ if and only if $p,q$ satisfy
\begin{enumerate}
\item $1/q \leq 1/p - \frac{k}{k+2}(1-\lam)$
\item  $1/q<\lam, 1/p>1-\lam.$
\end{enumerate}
Specifically, for such $p,q$,
\[ ||J_{Q_1,Q_2,\lam} f||_{\ell^q(\Z^{k+1})} \leq A_{p,q} ||f||_{\ell^p(\Z^{k+1})},\]
where  the constant $ A_{p,q}$ depends on $k, Q_1, Q_2,p,q,\lam$. 
\end{thm}

For the remaining values of $\lam$ in dimensions 3 and higher, interpolation yields a partial result:
\begin{cor}\label{Cor}
For $k \geq 3$ and $0< \lam \leq \lam_k$, $J_{Q_1,Q_2,\lam}$ extends to a bounded operator from  $\ell^p(\Z^{k+1})$ to $\ell^q(\Z^{k+1})$ if $p,q$ satisfy
\begin{enumerate}
\item $1/q < 1/p - (1- (\frac{k+2}{k})\lam)$
\item  $1/q<\lam, 1/p>1-\lam.$
\end{enumerate}
\end{cor}

Theorem \ref{dis_JQ_thm} is sharp in dimensions $k=1,2$ for all $0<\lam<1$. For dimensions $k \geq 3$, Theorem \ref{dis_JQ_thm} is sharp in the range $\lam_k < \lam < 1$, but for $0< \lam \leq \lam_k,$ the condition (i) in Corollary \ref{Cor} is weaker than the condition predicted by the continuous analogue, which is (i) in Theorem \ref{dis_JQ_thm}. 

In the case of $k=1$ with $Q_1(x) = Q_2(x) = x^2$, the result of Theorem \ref{dis_JQ_thm} follows from work of Stein and Wainger \cite{SW2}, \cite{SW3}, Oberlin \cite{Obe}, and Ionescu and Wainger \cite{IW}. In particular, the methods presented here in the higher dimensional setting are inspired by the beautiful ideas underlying the results of \cite{SW3}. In general, note that any two positive definite forms $Q,\tilde{Q}$ in $k$ variables are comparable, in the sense that there exist constants $c_1,c_2>0$ such that
$c_1Q(x) \leq \tilde{Q}(x) \leq c_2Q(x).$
As a result, the $(\ell^p, \ell^q)$  operator norms of $J_{Q_1,Q,\lam}$ and $J_{Q_1,\tilde{Q},\lam}$ are equivalent (up to constants); hence it would suffice to state Theorem \ref{dis_JQ_thm} in the case where the quadratic form in the denominator is the familiar form $|\cdot|^2$ given by the Euclidean norm.  However, the ability to choose the quadratic form in the denominator freely will play an important role in the method of proof; in keeping with this flexibility, we state the theorem in greatest generality. 
Finally, we note that $J_{Q_1,Q_2,\lam}$ may be seen as a translation invariant model of a discrete analogue of fractional integration on the Heisenberg group; the author will address the latter quasi-translation invariant operator in a later paper. 

It is illustrative to compare Theorem \ref{dis_JQ_thm} to the $(L^p, L^q)$ behavior of the continuous analogue of $J_{Q_1,Q_2,\lam}$. 
Even in the Euclidean setting, fractional integral operators over submanifolds are not yet completely understood; in particular it is not known precisely how the $(L^p, L^q)$ inequalities of such operators depend on the specific geometry of the underlying submanifold. However, one special case that is understood is that of fractional integration along a paraboloid in $\R^k$, given by
\beq\label{cts_J_dfn}
\mathcal{J}_{\lam}f(x,t) = \int_{\R^k} f(x-y,t-|y|^2)\frac{dy}{|y|^{k\lam}},
\eeq
for $0< \lam < 1$.
In the case $k=1$, Christ \cite{Chr88} proved sharp results for $\mathcal{J}_\lam$ via a Littlewood-Paley decomposition. His methods may be extended to higher dimensions and more general positive definite quadratic forms in $k$ variables, resulting in the following:
let $\Gamma$ be the closed triangle in $[0,1]^2$ defined by the vertices $(0,0)$, $(1,1)$, $(1/p_k, 1/q_k)$ where $1/p_k = (k+1)/(k+2)$, $1/q_k = 1/(k+2)$. 
Then for $0< \lam < 1$, $\mathcal{J}_{\lam}$ is a bounded operator from $L^p(\R^{k+1})$ to $L^q(\R^{k+1})$ precisely for those pairs $p,q$ such that $(1/p,1/q) \in \Gamma$ and $1/q = 1/p-\frac{k}{k+2}(1-\lam).$

\subsection{Outline of the proof}
The proof of Theorem \ref{dis_JQ_thm} uses techniques from the circle method of Hardy and Littlewood, theta functions, and exponential sums, as well as more classical analytic techniques of complex interpolation, Fourier transform methods, and oscillatory integral estimates. Because of the intricacy of the method, we outline  the main ideas of the proof here, before turning to a technical presentation.

Proof that conditions (i) and (ii) in Theorem \ref{dis_JQ_thm} are necessary for the $(\ell^p, \ell^q)$ boundedness of the operator follows from facts about the average number of representations of a positive integer by a given positive definite quadratic form; see the appendix in Section \ref{nec_sec}.
The main difficulty lies in proving the sufficiency of conditions (i) and (ii).

By choosing the quadratic form $Q_2$ to equal $Q_1,$ we may immediately reduce consideration to an operator of the form
\beq\label{dis_JQ_dfn_Q} 
 J_\lam f(n,t) = \sum_{\bstack{m \in \Z^k}{m \neq 0}} \frac{f(n-m,t-Q(m))}{Q(m)^{k\lam/2}}.
 \eeq
Since it is translation invariant, $J_\lam$ is a Fourier multiplier operator,
\beq\label{J_F_mult}
 (J_\lam f)\hat{\:} (\theta,\phi) = m_\lam(\theta,\phi) \hat{f}(\theta,\phi),
 \eeq
where the multiplier $m_\lam(\theta,\phi)$ is given for $\theta \in [0,1]$, $\phi \in [0,1]^k$ by
\beq\label{mult_dfn0}
 m_\lam (\theta, \phi) = \sum_{\bstack{m \in \Z^k}{m \neq 0}} \frac{e^{-2\pi i (Q(m) \theta + m \cdot \phi)}}{Q(m)^{k\lam/2}}.
 \eeq
(Here we abuse notation slightly by writing the spectral variables in the order $(\theta,\phi)$, where $\theta$ corresponds to $t$ and $\phi$ to $n$; this is motivated by the privileged role $\theta$ plays, due to the quadratic nature of its coefficient in the phase of $m_\lam$.)
We will show that the multiplier $m_\lam$ may be approximated (up to acceptable error) by the multiplier 
\beq\label{nu_int}
 \nu_{\lam}(\theta,\phi) = \int_0^1 \Theta_Q(y+i\theta,\phi) y^{k\lam/2-1} dy ,
 \eeq
where
\beq\label{theta_dfn0}
\Theta_Q(z,\phi) = \sum_{m \in \Z^k}e^{-2\pi Q(m)z}e^{-2\pi i m \cdot \phi}
\eeq
is a twisted theta function defined for $z\in \C, \phi \in [0,1]^k$ and absolutely convergent for $\Re(z)>0$. Here we have used the fact that we chose $Q_2=Q_1 =Q$ in order to write $m_\lam$ in the form $\nu_\lam$.

We proceed by decomposing the multiplier $\nu_\lam$ in two ways, first as a sum over dyadic portions of the integral (\ref{nu_int}), and second, in terms of the Diophantine properties of the spectral variables $\theta$ and $\phi$, according to ``major'' and ``minor'' arcs. Intuitively, when $\theta$ may be closely approximated by a rational number $a/q$ where $1 \leq a \leq q$, $(a,q)=1$, and $q$ is ``small,'' $\theta$ is said to lie in a major arc, and otherwise $\theta$ lies in a minor arc. In this way $J_\lam$ may be decomposed on the spectral side into an infinite family of ``major'' and ``minor'' operators, each with certain useful arithmetic properties.

Ultimately, after a careful parameterization of the major arcs, this results in a representation (up to acceptable error) of the contribution of the major arcs to the operator $J_\lam$ as an infinite sum of a family of operators \beq\label{J_decomp}
\sum_{s=0}^\infty \mathcal{B}_{\lam,s}, 
\eeq
in which each component is defined by a corresponding Fourier multiplier $B_\lam(s; \theta,\phi)$,
\[ ( \mathcal{B}_{\lam,s} f)\hat{\;} (\theta, \phi) = B_\lam(s;\theta, \phi) \hat{f}(\theta,\phi).\]
(A similar decomposition holds for the contribution of the minor arcs to $J_\lam,$ and a procedure similar to that described below is followed.) 
It is then sufficient to show for each index $s \geq 0$ that for all appropriate $\lam, p,q$,
\beq\label{B_lam_s}
 || \mathcal{B}_{\lam,s} f||_{\ell^{q}} \leq C2^{-\del(\lam)s} ||f||_{\ell^p} 
 \eeq
for some $\del(\lam)>0$, resulting in a finite $(\ell^p, \ell^q)$ norm for the sum (\ref{J_decomp}).

Our approach to proving (\ref{B_lam_s}) is by complex interpolation of the operators $\mathcal{B}_{\lam,s}$ in the region $-2/k \leq \Re(\lam) \leq 1$; our argument elaborates upon the basic principles underlying the results of \cite{SW3}. Specifically, we prove bounds of the following form: there exist nontrivial constants $\al, \be>0$ such that for $\lam = 1 + i\gamma$ with $\ga \neq 0$,
\beq\label{22}
 || \mathcal{B}_{\lam,s} f||_{\ell^{2}} \leq A2^{-\al s} ||f||_{\ell^2}, 
 \eeq
and for $\lam = -2/k + i\gamma$,
\beq\label{1inf}
 || \mathcal{B}_{\lam,s} f||_{\ell^{\infty}} \leq B2^{+ \be s} ||f||_{\ell^1} .
 \eeq
If the growth due to $2^{+\be s}$ is sufficiently small with respect to the decay provided by $2^{-\al s}$, complex interpolation results in a bound of the form (\ref{B_lam_s}) for $\lam$ in a sufficiently large range to prove Theorem \ref{dis_JQ_thm}.

The $(\ell^2,\ell^2)$ bounds of the form (\ref{22}) on the line $\Re(\lam)=1$ follow from $\ell^\infty$ bounds for the corresponding multipliers $B_\lam(s;\theta,\phi)$. On the line $\Re(\lam)=-2/k$, the original operator $J_\lam$ is not easily interpreted, and certainly does not exhibit bounded behavior; however, each individual piece $\Bcal_{\lam,s}$ within the decomposition may be bounded from $\ell^1$ to $\ell^\infty$, albeit with visibly exponential growth in the norm. Interestingly, these $(\ell^1,\ell^\infty)$ bounds follow from estimates not for the size of the multipliers $B_\lam(s;\theta,\phi)$ but of their Fourier coefficients; these estimates in particular must be handled with great care. 

Along the way, complications arise due to the singularity of certain oscillatory integrals, forcing us to adopt a substantially more intricate decomposition  than indicated in the simplistic representation (\ref{J_decomp}), but this outlines the main framework of the proof. 
Throughout, the beautiful structure of the theta function $\Theta_Q(y+i\theta,\phi)$ shines through.
In particular, the theta function provides a transformation law (effectively the Jacobi inversion formula) that allows us to separate each  multiplier $B_\lam(s;\theta,\phi)$ into an arithmetic part and an integral, depending explicitly on the approximation of $\theta$ by a rational number $a/q$ and of $\phi$ by a $k$-tuple of rational numbers $(b_1/q, \ldots, b_k/q)$. The aim, on both the line $\Re(\lam)=1$ and the line $\Re(\lam)=-2/k$, is to extract optimal cancellation from both the arithmetic part (an exponential sum) and the integral. 

Discrete fractional integral operators of the form $J_\lam$ defined in (\ref{dis_JQ_dfn_Q}) and the more general form $I_\lam^P$ defined in (\ref{IP_lam}) are closely related to translation invariant discrete singular Radon transforms of critical degree defined by
\beq\label{TP_dfn}
T_Pf(n) = \sum_{\bstack{m \in \Z^{k_1}}{m \neq 0}} f(n-P(m))K(m),
\eeq
where $K$ is a Calder\'{o}n-Zygmund kernel on $\R^{k_1}$ and $P=(P_1, \ldots, P_{k_2})$ is a polynomial mapping from $\Z^{k_1}$ to $\Z^{k_2}$.
A recent deep result of Ionescu and Wainger \cite{IW} proves that $T_P$ is bounded on $\ell^p(\Z^{k_2})$ for all $1<p< \infty$, with operator norm dependent only on $p$, the dimension $k_1$, and the degree of the polynomial $P$. This result will play a role in the proof of Theorem \ref{dis_JQ_thm}, but it also highlights a contrast between $(\ell^p, \ell^q)$ results for discrete fractional Radon transforms, such as $J_\lam$, and $(\ell^p, \ell^p)$ results for Radon transforms with singularities of critical degree, such as $T_P$. 

Indeed, it might appear at first glance that $\ell^2$ bounds of the form (\ref{22}) should already be known from the work of Ionescu and Wainger on $T_P$. Superficially this is true---however, in the case of $T_P$, it was sufficient to prove bounds of the form (\ref{22}) with arbitrarily small decay exponent $\al$ with respect to $s$. In contrast, in our application we must extract maximal decay in $s$ in order to counteract the growth present in the accompanying bounds (\ref{1inf}). This explains the care with which we must treat each and every estimate that arises, which might otherwise seem unnecessarily parsimonious. Furthermore, this is the reason that the result of Ionescu and Wainger holds for operators $T_P$ with a polynomial map $P$ of any degree, while Theorem \ref{dis_JQ_thm} considers only quadratic forms: due to the demands of the interpolation argument, 
we require optimal (square-root) cancellation in certain exponential sums, which in this setting is only available in the quadratic case. 

\subsection{Outline of the paper}
In Section \ref{sec_Prelim_FR} we fix notation and prove several number theoretic lemmas on  exponential sums and theta functions associated to the quadratic form $Q$. In Section \ref{sec_decompJ} we decompose the operator $J_\lam$ dyadically and according to major and minor arcs. However, singularities of certain integrals cause us to adopt a more complicated double decomposition of the major arc operators, which we motivate in Section \ref{sec_dd_motiv}. We then carry out this double decomposition and proceed to estimate the resulting multipliers in Sections \ref{sec_double_decomp}, \ref{sec_maj_rem} and \ref{sec_min_full}. In an appendix in Section \ref{nec_sec} we record examples that show that conditions (i) and (ii) of Theorem \ref{dis_JQ_thm} are necessary for the $(\ell^p, \ell^q)$ boundedness of the operator.

\subsection*{Note} This is a corrected version of the original paper, which appeared in {\emph{Duke Math. Journal}} {\bf{161}} no. 1 (2012) 69-106, and contained an error in equation (83). The statements of Propositions 3, 6, 7, and Theorem 1 have been corrected, and Corollary 1.1 has been added. (September 2012)

\section{Preliminaries}\label{sec_Prelim_FR}
\subsection{Notation}
The discrete Fourier transform of a function $f \in \ell^1(\Z^k)$ is defined by 
\[ \hat{f}(\theta) = \sum_{n \in \Z^k} f(n)e^{-2\pi i n \cdot \theta}\]
where $\theta \in [0,1]^k$.
The Fourier inverse of a periodic function $h(\theta)\in L^2_{loc}(\R^k)$ is defined by 
\[ \check{h}(n) = \int_{[0,1]^k} h(\theta) e^{2\pi i \theta \cdot n } d \theta,\]
for any $n \in \Z^k$.  
If $f(n_1,n_2)$ acts on $\Z^{k_1}\times \Z^{k_2}$, we will denote by $f\hat{^{n_2}}(n_1,\theta)$ the ``partial'' Fourier transform with respect to the variable $n_2$ alone, where $n_1 \in \Z^{k_1}$ and $\theta \in [0,1]^{k_2}$.

Given a positive definite quadratic form $Q$ in $k$ variables with integer coefficients, we may write $Q(x) = \frac{1}{2}A[x] = \frac{1}{2} x^t A x,$
where $A$ is a real, positive definite, $k \times k$ symmetric matrix with integer entries and even diagonal entries.  
We assume $Q$ acts nontrivially in all $k$ variables, i.e. is of rank $k$, so that $A$ is invertible; moreover, there exists an orthogonal matrix $P$ and a diagonal matrix $D$ such that $A =P^t D P$. 
We define the adjoint quadratic form to $Q$ by $Q^*(x) =\frac{1}{2} A^{-1}[x]$. 
In number theoretic terms, the quadratic form $|x|^2$ is quite special, but it arises in analytic problems as the most familiar case, and thus we will refer to it in this paper as the {\emph{generic quadratic form}}.
We will also use $A$ to denote a constant (possibly dependent on the dimension $k$), which may change from one occurrence to the next; it will be clear from context whether $A$ denotes the matrix or a constant.

\subsection{Gauss sums}
An important aspect of the method we present involves extracting arithmetic information from Fourier multipliers in the form of exponential sums, and then harvesting cancellation from these sums.
In particular, we will require the Gauss sum defined in terms of a positive definite quadratic form $Q$ by
\beq\label{Gauss_sum_dfnQ}
 S_Q(a,b;q) =  \sum_{r \modd{q}^k} e^{-2\pi i Q(r) a/q}e^{-2\pi i r \cdot b/q},
 \eeq
 where the notation $r \modd{q}^k$ (or equivalently $r \in (\Z/q\Z)^k$) indicates the collection of $k$-tuples $r = (r_1, \ldots, r_k)$ with each $1 \leq r_j \leq q$. Similarly, $r \con a \modd{q}^k$ with $r,a \in \Z^k$ indicates that for each $1 \leq j \leq k$, $r_j \con a_j \modd{q}$. 
The classical bounds we require are as follows (see also Lemma 20.12 of \cite{IK}):
\begin{lemma}\label{Gauss_sum_lemmaQ}
For any $(a,q)=1$, $b \in \Z^k$, the Gauss sum $S_Q(a,b;q)$ is bounded by
\beq\label{Gauss_boundQ}
|S_Q(a,b;q)| \leq c_Q q^{k/2},
\eeq
where the constant $c_Q$ depends only on $Q$ and $k$.
 Furthermore, for any $l_1 \in \Z, l_2 \in \Z^k$,
\beq\label{Gauss_avg_boundQ}
 \left| \sum_{\bstack{a=1}{(a,q)=1}}^q \sum_{b \modd{q}^k} S_Q(a,b;q) e^{-2\pi i l_1 a/q}e^{-2\pi i l_2 \cdot b/q} \right| \leq q^{k+1}.
 \eeq
\end{lemma}
To prove (\ref{Gauss_boundQ}), we square the sum, noting that if $A$ is the matrix representing $Q$ then:
\begin{eqnarray}
|S_Q(a,b;q)|^2 &=& \sum_{r,s \modd{q}^k} e(-\frac{a}{q} (Q(r) - Q(s))) e(-(r-s)\cdot \frac{b}{q}) \nonumber \\
& = & \sum_{u,s \modd{q}^k} e(-\frac{a}{q} (Q(u)+s^tAu))e(- u \cdot \frac{b}{q}) \nonumber \\
& = & \sum_{u \modd{q}^k} e(-\frac{1}{q}(aQ(u) + u \cdot b)) \sum_{s \modd{q}^k} e(-\frac{a}{q} (s^t Au)) \nonumber \\
& \leq & q^k | \{ u \modd{q}^k : Au \con 0 \modd{q}^k \} |.\label{Au}
\end{eqnarray}
It remains to count the number of $ u \modd{q}^k$ such that $v := Au \con 0 \modd{q}^k$. Let $\del = \det{A}$. Since $v \con 0 \modd{q}^k$ then $\del u \con 0 \modd{q}^k$; hence $\del u/q$ is an integer vector. Set $d=(\del,q)$. Then it remains true that $(\frac{\del/d}{q/d})u = \del u/q$ is an integer vector, but now $(\del/d, q/d)=1$; from this it follows that $q/d$ must divide each coordinate of $u$, so that $u \con 0 \modd{q/d}^k$. Recall that by definition $u$ runs modulo $q$ in each coordinate, so if each coordinate must be divisible by the integer $q/d$, there are at most $d$ choices for each coordinate, and hence $d^k$ choices for $u$ in total. But recall furthermore that $d | \del$, and hence $d \leq \del$, so there are at most $\del^k$ choices of $u$ such that $Au \con 0 \modd{q}^k$.
Thus in (\ref{Au}) we obtain $|S_Q(a,b;q)|^2 \leq \del^k q^k$. Taking square roots, this yields a bound of size $|\det A|^{k/2}q^{k/2}$ for $S_Q$, proving (\ref{Gauss_boundQ}).

For the second bound (\ref{Gauss_avg_boundQ}), we first expand the Gauss sum, obtaining
\beq\label{expand}
\sum_{\bstack{a=1}{(a,q)=1}}^q  \sum_{r \modd{q}^k}  e^{-2\pi i Q(r) a/q} e^{-2\pi i l_1 a/q}\sum_{b \modd{q}^k}e^{-2\pi i r \cdot b/q} e^{-2\pi i l_2 \cdot b/q}.
\eeq
The innermost sum over $b$ contributes $ q^k$ if $r \con - l_2 \modd{q}^k$ and vanishes otherwise, so that (\ref{expand}) equals
\[ q^k \sum_{\bstack{a=1}{(a,q)=1}}^q  e^{-2\pi i (Q(-l_2)+l_1)a/q}.\]
This is trivially bounded by $q^{k+1}$ in absolute value, and this completes the proof.

\subsection{Theta function of a quadratic form}
Recall the twisted theta function associated to the positive definite quadratic form $Q$, defined by 
\beq\label{theta_phi_Q_dfn}
 \Theta_Q(z,\phi) = \sum_{m \in \Z^k}e^{-2\pi Q(m)z}e^{-2\pi i m \cdot \phi},
\eeq
where $z \in \C$, $\phi \in [0,1]^k$. 
This is absolutely convergent for $\Re(z) >0$ and uniformly convergent in closed half-planes $\Re(z) \geq \del >0$.
The key identity we require for $\Theta_Q(y+i\theta,\phi)$ is essentially the Jacobi inversion formula, specialized according to the approximation of $\theta,\phi$ by rational numbers.
\begin{lemma}\label{G_r_transf_lawQ}
 For $y>0$, $\theta \in [0,1], \phi \in [0,1]^k$ with $\theta =a/q +\al$, $\phi =b/q + \be$, 
 \[ \Theta_Q(y+i\theta,\phi)  = \frac{1}{q^k |A|^{1/2}(y+i\al)^{k/2}} \sum_{m \in \Z^k} S_Q(a, b-m;q) e^{-2\pi Q^*(m/q + \be)/(y+i\al)}.\]
\end{lemma}
To derive this, write $m=lq + r$ in the definition of $\Theta_Q(z,\phi)$, where $l \in \Z^k, r \in (\Z/q\Z)^k$, so that
\begin{eqnarray}
\Theta_Q(y+i\theta,\phi)
	&=& \sum_{l \in \Z^k} \sum_{r \in (\Z/q\Z)^k} e^{-2\pi Q(lq+r)(y+i\al + ia/q)} e^{-2\pi i (lq+r) \cdot (b/q + \be)} \nonumber \\
	& = & \sum_{r \in (\Z/q\Z)^k} e^{-2\pi i Q(r)a/q} e^{-2\pi i r \cdot b/q} G_r(y+i\al,\be),\label{theta1}
\end{eqnarray}
where
\[G_r(z,\beta) = \sum_{l \in \Z^k} e^{-2\pi Q(lq +r)z}e^{-2\pi i (lq+r) \cdot \beta}.\]
Applying Poisson summation to the function pair 
\begin{eqnarray*}
 h(x) &=& e^{-2\pi Q(xq +r)y}e^{-2\pi i (xq+r) \cdot \beta} ,\\
\hat{h}(\xi) &=& \frac{e^{2\pi i r \cdot \xi/q}}{q^{k}|A|^{1/2}y^{k/2}} e^{-2\pi Q^*(\frac{\xi + \be q}{q\sqrt{y}})},
\end{eqnarray*}
one obtains for real $z=y>0$ that
\beq\label{G_eqn}
G_r(z,\beta) 
	= \frac{1}{q^k|A|^{1/2}z^{k/2}}\sum_{l \in \Z^k} e^{2\pi i r \cdot l/q}e^{-2\pi Q^*(l/q +\be)/z},
	\eeq
	where $Q^*$ is the adjoint quadratic form to $Q$. This identity continues to hold for all $\Re(z)>0$ by analytic continuation, using the principal branch of the logarithm to define $z^{-k/2}$.
Recalling the definition (\ref{Gauss_sum_dfnQ}) of the Gauss sum $S_Q$, the lemma then follows immediately from (\ref{theta1}) and (\ref{G_eqn}).

\section{A basic decomposition for $J_{\lam}$}\label{sec_decompJ}
\subsection{Main proposition}
We now turn to the proof proper of Theorem \ref{dis_JQ_thm}. In order to prove the theorem it is sufficient to prove the following proposition.
\begin{prop}[Main Proposition]\label{key_J_prop}
For $\max(\frac{2}{k+4}, \frac{k}{2k+2})<\lam < 1,$ and for all $1/p + 1/p' =1$ with $1/p'=1/p - k(1-\lam)/(k+2)$,
\[ ||J_\lam f||_{\ell^{p'}} \leq A_p ||f||_{\ell^p}.\]
\end{prop}
That this is sufficient may be seen by standard interpolation arguments, which we outline briefly. Recall the  result of Ionescu and Wainger \cite{IW} that  discrete singular Radon transforms of the form (\ref{TP_dfn}) are bounded operators on $\ell^p$ for $1< p < \infty$. When $\Re(\lam)=1$ with $\Im(\lam) \neq 0$, $J_\lam$ is of this form, and hence is bounded on $\ell^p$ for all $1 < p < \infty$. Complex interpolation of this result with Proposition \ref{key_J_prop} for the family of operators $e^{\lam^2-1}(1-\lam) J_\lam$ removes the conjugacy condition $1/p+1/p'=1$ for  $\lam > \max(\frac{2}{k+4}, \frac{k}{2k+2})$. For $k=1,2$ the allowed range is $2/(k+4) < \lam <1$. The value $\lam=2/(k+4)$ is the ``crossover point'' at which condition (ii) in Theorem \ref{dis_JQ_thm} becomes more restrictive than condition (i). For $k=1,2$, to complete the proof of Theorem \ref{dis_JQ_thm} in the region $0< \lam \leq 2/(k+4)$, it suffices to obtain all pairs $(1/p, 1/q)$ satisfying condition (ii). We may do this by interpolating the result of Proposition \ref{key_J_prop} for $\lam = \frac{2}{k+4} + \ep$ with the trivial $\ell^1 \maps \ell^\infty$ bound for $J_\lam$ on the line $\Re(\lam)=0$, and then using the inclusion principle for $\ell^p$ spaces, namely $\ell^{q_2} \subset \ell^{q_1}$ if $q_2 < q_1$, and finally taking adjoints. For $k \geq 3$, to obtain Corollary \ref{Cor}, we follow the same procedure but interpolate the result of Proposition \ref{key_J_prop} for $\lam = \frac{k}{2k+2} + \ep$ with the trivial $\ell^1 \maps \ell^\infty$ bound for $J_\lam$ on the line $\Re(\lam)=0$. In this respect we note that $1-(\frac{k+2}{k})\lam$ is the linear function in $\lam$ that equals $\frac{k}{k+2}(1-\lam)$ when $\lam = \frac{k}{2k+2}$ and equals $1$ when $\lam =0$. 
Note that we could prove the result of Theorem \ref{dis_JQ_thm} for all $0<\lam<1$ and all $k \geq 3$ if we could extend the validity of Proposition \ref{key_J_prop} to $\frac{2}{k+4} < \lam < 1$ for all $k \geq 3$.

\subsection{Dyadic decomposition of the multiplier}\label{sec_dyadic_decomp}
To prove the main proposition, we begin by reducing the multiplier of $J_\lam$ to the integral form (\ref{nu_int}), which we then decompose dyadically. Recall from (\ref{mult_dfn0}) that $J_\lam$ has Fourier multiplier
\[ m_\lam (\theta, \phi) = \sum_{\bstack{m \in \Z^k}{m \neq 0}} \frac{e^{-2\pi i (Q(m) \theta + m \cdot \phi)}}{Q(m)^{k\lam/2}}.\]
The identity
\[
 Q(m)^{-k\lam/2} = c_{k,\lam} \int_0^\infty e^{-2\pi Q(m) y}y^{k\lam/2-1}dy,\]
 with $c_{k,\lam} = (2\pi)^{k\lam/2}/\Ga(k\lam/2)$,
 allows us to write
\[ m_\lam (\theta,\phi) = c_{k,\lam} \int_0^\infty \Theta_0(y+i\theta, \phi) y^{k\lam /2-1}   dy,\]
where  $\Theta_0(z,\phi) = \Theta_Q(z,\phi)-1$, and $\Theta_Q$ is the theta function (\ref{theta_phi_Q_dfn}).
The portion of the integral over $[1,\infty)$ has absolutely convergent Fourier series and hence the operator with multiplier corresponding to this integral is bounded from $\ell^p$ to $\ell^q$ for $1 \leq p \leq q$. Furthermore, in the portion of the integral over $[0,1]$, introducing the $m=0$ term of the theta function contributes only an $O(1)$ error term to the theta function, which will be subsumed by the error $E_\lam$ encountered later (see equation (\ref{remainder_dfnQ})).
Thus it is sufficient to consider the multiplier
\beq\label{nu_dfn}
 \nu_{\lam}(\theta,\phi) = \int_0^1 \Theta_Q(y+i\theta,\phi) y^{k\lam/2-1} dy ,
 \eeq
which we further decompose  as $\nu_{\lam} =\sum_{j=1}^\infty \nu_{\lam,j}$, where
 \beq\label{nu_j_dfn}
 \nu_{\lam,j}(\theta, \phi) = \int_{2^{-j}}^{2^{-j+1}}  \Theta_Q(y+i\theta,\phi) y^{k\lam/2-1} dy .
 \eeq

\subsection{Major and minor arcs}\label{sec_mm_decomp}
Our goal is to decompose the multiplier according to the Diophantine properties of the spectral variables $(\theta,\phi)$. We will do so by decomposing the unit cube $[0,1] \times [0,1]^k$ into major and minor arcs, following the nomenclature of the circle method.
By Dirichlet's approximation principle, given any $N \geq 1$, for every $\theta \in [0,1]$ there exist integers $a,q$ with $1 \leq q \leq N$, $1 \leq a \leq q$ and $(a,q)=1$, such that 
\[ | \theta - a/q| \leq \frac{1}{ qN}.\]
For each $j$ indexing a dyadic range of the multiplier $\nu_{\lam,j}(\theta,\phi)$, we will apply this principle (with $N=2^{j/2}$) to choose a rational approximation to $\theta$. Moreover, we will define the major arcs (for that fixed $j$) to correspond to those $\theta$ approximated by a rational $a/q$ with $q$ small (specifically $1 \leq q \leq 2^{j/2-10})$ and the minor arcs to be the complement of the collection of major arcs (for that fixed $j$). In fact, it will be convenient to use two parameters: $j$, indexing the dyadic decomposition of the multiplier integral, and $s$, indexing  dyadic ranges of denominators $q$. 

Precisely, 
we define the major arcs for each $j$ as follows: for $0 \leq s \leq j/2-10$, $2^s \leq q< 2^{s+1}$, $1\leq a \leq q$ with $(a,q)=1$, and $b=(b_1,\ldots, b_k)$ with $1 \leq b_i \leq q$, let 
\beq\label{Maj}
 \M_{j,s}(a,b;q) = \{(\theta, \phi): |\theta - a/q| \leq \frac{1}{2^s2^{j/2}}, |\phi_i - b_i/q| \leq \frac{1}{2q}, i=1,\ldots, k \}.
 \eeq
 Later we will let $\M_j(a,b;q)$ stand for the union over $s \leq j/2-10$ of $\M_{j,s}(a,b;q)$.
The minor arcs for each $j,s$ are defined to be the complement of the union of the collection of major arcs (for that fixed $j,s$ pair) in $[0,1]^{k+1}$. We will always use the notation $\theta= a/q +\al$ and $\phi = b/q + \beta$, where $b/q$ represents $(b_1/q, \ldots, b_k/q)$ and $\al\in [0,1]$, $\be \in [0,1]^k$.

We note that for fixed $j$ but varying $s$, the major arcs are disjoint. For suppose $j$ is fixed and consider $\M_{j,s_1}(a_1,q_1)$ and $\M_{j,s_2}(a_2,q_2)$, where $s_1 \neq s_2$. Then if these two arcs intersected, they would intersect in all coordinates, and thus in particular in the first coordinate, and we would have 
\[ \frac{1}{4 \cdot 2^{s_1} 2^{s_2}} < \frac{1}{q_1 q_2} \leq \left| \frac{a_1}{q_1} - \frac{a_2}{q_2}\right| \leq \left( \frac{1}{2^{s_1}2^{j/2}} + \frac{1}{2^{s_2}2^{j/2}} \right) \leq \frac{2}{2^{j/2}}\max \left( \frac{1}{2^{s_1}}, \frac{1}{2^{s_2}} \right),\]
which is impossible since both $s_1, s_2 \leq j/2-10$.

On the other hand, later we will switch the order of summation, letting $s$ run to infinity and summing over $j \geq 2s+20$. In this instance we will need to know that for fixed $s$ but varying $j$ the major arcs corresponding to $a_1/q_1 \neq a_2/q_2$ are disjoint. This is proved similarly: if $s$ is fixed and $\M_{j_1,s}(a_1,q_1)$ and $\M_{j_2,s}(a_2,q_2)$ are given, where $(a_1,q_1) \neq (a_2,q_2)$, then if the arcs intersected we would have
 \[ \frac{1}{4 \cdot 2^{2s}} < \frac{1}{q_1 q_2} \leq \left|\frac{a_1}{q_1} - \frac{a_2}{q_2} \right| \leq \left( \frac{1}{2^s2^{j_1/2}} + \frac{1}{2^s 2^{j_2/2}} \right) \leq \frac{2}{2^{s}}\max \left( \frac{1}{2^{j_1/2}}, \frac{1}{2^{j_2/2}} \right).\]
But both $j_1, j_2 \geq 2s +20$, so this is impossible.

\subsection{The approximate identity}
We will now derive an approximate identity for the theta function $\Theta_Q(z,\phi)$ that will allow us to harvest the arithmetic information encoded in the rational approximations to $\theta,\phi$; this identity holds for $(\theta,\phi)$ in either a major or a minor arc.
\begin{prop}[Approximate identity]\label{approx_id_no_derivQ}
Suppose $1 \leq q \leq 2^{j/2}$, $1 \leq a \leq q$, $(a,q)=1$, and $1 \leq b_i \leq q$ for $i=1,\ldots, k$. Then for $(\theta,\phi)$ and $a,b,q$ such that $|\theta-a/q| \leq 1/(q2^{j/2})$ and $|\phi_i - b_i/q| \leq 1/(2q)$ for $i=1, \ldots, k$, we have:
\beq\label{approx_id_eqn_JQ}
\Theta_Q(y + i\theta,\phi) =  \frac{S_Q(a,b;q)}{q^k|A|^{1/2}} \cdot \frac{e^{-2\pi Q^*(\be)/(y+i\al)}}{(y+i\al)^{\frac{k}{2}}} + E_{\lam}(y+i\theta,\phi) ,
\eeq
where the remainder term is given by 
\beq\label{remainder_dfnQ}
 E_{\lam}(y+i\theta,\phi) = \sum_{\bstack{m \in \Z^k}{m \neq 0}} \frac{S_Q(a,b-m;q)}{q^k|A|^{1/2}} \cdot \frac{e^{-2\pi Q^*(m/q+\be)/(y+i\al)}}{(y+i\al)^{\frac{k}{2}}}.
 \eeq
Furthermore, if $2^{-j} \leq y \leq 2^{-j+1}$ then
\beq\label{E_bdQ}
 |E_{\lam}(y+i\theta,\phi)|  = O(y^{-k/4}).
 \eeq
\end{prop}

Recall from Lemma \ref{G_r_transf_lawQ} that $\Theta_Q$ satisfies the transformation law
\[ \Theta_Q(y+i\theta,\phi)  = \frac{1}{q^k |A|^{1/2}(y+i\al)^{k/2}} \sum_{m \in \Z^k} S(a, b-m;q) e^{-2\pi Q^*(m/q + \be)/(y+i\al)}.\]
The main term of (\ref{approx_id_eqn_JQ}) is simply the $m=0$ term; $E_{\lam}(y+i \theta,\phi)$ comprises the remainder of the sum. Note that
by the Gauss sum bound (\ref{Gauss_boundQ}) of Lemma \ref{Gauss_sum_lemmaQ},
\beq\label{E_orig}
 |E_{\lam}(y+i\theta,\phi)| \leq \frac{C}{q^{k/2}|y+i\al|^{k/2}} | \sum_{\bstack{m \in \Z^k}{m \neq 0}} e^{-2\pi Q^*(m/q+\be)/(y+i\al)}|.
\eeq
By a simple comparison we may replace the quadratic form $Q^*$ by the generic quadratic form $|\cdot|^2$, so that  the sum over $m \neq 0$ is controlled by a product of sums, 
\beq\label{prods}
 \prod_{i=1}^k |\sum_{m_i } e^{-c|m_i/q+ \be_i|^2/(y+i\al)} |
\leq \prod_{i=1}^k \sum_{m_i  } e^{-c|m_i+\be_i q|^2y/(q^2(y^2+\al^2))},
\eeq
where in each case not all $m_i$ are simultaneously zero.
Recall that in either a major or minor arc, $|\al| \leq q^{-1}2^{-j/2}$ and $|\be_i| \leq c(2q)^{-1}$, where $c$ depends only on the dimension. Thus under the hypotheses of Proposition \ref{approx_id_no_derivQ}, $q \leq c y^{-1/2}$ and $q|\al| \leq cy^{1/2}$ and thus setting $u  = y/(q^2(y^2 + \al^2))$,  it follows that $u \geq c \cdot 1$. Since in all cases at least one factor in (\ref{prods}) omits $m_i=0$, it follows that (\ref{prods}) is $O(u^{-k/4})$. We record this for future reference as 
\beq\label{general_m_sum_bound}
 |\sum_{m\neq 0}  e^{-2\pi Q^*(m/q + \be)/(y+i\al)}|  \leq A y^{-k/4} q^{k/2}(y^2 + \al^2)^{k/4}.
\eeq 
(Note that this continues to hold under the slightly more general condition that $|\be_i| \leq 3/(4q)$, as we will assume later.)
In conclusion, 
\[ |E_\lam(y+i \theta,\phi)|  = O(y^{-k/4}),\]
and this completes the proof of Proposition \ref{approx_id_no_derivQ}.

\subsection{A major/minor decomposition of the operator}\label{sec_J_lam_hard}
We now define the multipliers according to the major/minor decomposition of Section \ref{sec_mm_decomp}.
We begin by considering the multiplier $M_\lam$ corresponding to the main term in the approximate identity (\ref{approx_id_eqn_JQ}) for $\Theta_Q$. Let $\chi$ denote the characteristic function of $[-1,1]$ and $\uchi$  the characteristic function of $[-1,1]^k$. Using the notation $\theta = a/q + \al$, $\phi =b/q+\be$, the multiplier $M_\lam$ is given by:
\begin{eqnarray}
M_\lam(\theta,\phi) & = & \sum_{j=1}^\infty \sum_{s=0}^{j/2-10} \sum_{2^s \leq q < 2^{s+1}} \sum_{\bstack{a=1}{(a,q)=1}}^q  \sum_{b  \in (\Z/q\Z)^k}
\frac{S_Q(a,b;q)}{q^k|A|^{1/2}}\chi(2^{j/2}2^s(\theta -a/q)) \nonumber \\
	& & \cdot \; \uchi(2q(\phi-b/q))  \int_{2^{-j}}^{2^{-j+1}} \frac{y^{k\lam/2-1} e^{-2\pi Q^*(\be)/(y+i\al)}}{(y+i\al)^{\frac{k}{2}}} dy  \label{M_mult_dfn} \\
	& = & \sum_{s=0}^\infty \sum_{j \geq 2s +20} \cdots \nonumber \\
	& = & \sum_{s=0}^\infty B_\lam (s; \theta, \phi), \nonumber
\end{eqnarray}
say.
Summing first in $j$ and using the support of the cut-off factors, we may write
\begin{multline}\label{B_shifted}
B_\lam (s; \theta, \phi)  =  \sum_{2^s \leq q < 2^{s+1}} \sum_{\bstack{a=1}{(a,q)=1}}^q  \sum_{b \in (\Z/q\Z)^k}
\frac{S_Q(a,b;q)}{q^k|A|^{1/2}}  \\
	\cdot \uchi(2q(\phi-b/q)) \int_{\rho(s,\al)}^{2^{-2s-19}} \frac{y^{k\lam/2-1} e^{-2\pi Q^*(\be)/(y+i\al)}}{(y+i\al)^{\frac{k}{2}}} dy,
\end{multline}
where $\rho(s,\al)$ represents $2^{-J}$ where $J$ is the largest $j \geq 2s +20$ such that $|\al| = |\theta -a/q| \leq 2^{-j/2}2^{-s}$, i.e. the largest $j$ such that $2^{-j} \geq 2^{2s}\al^2.$ Therefore
\beq\label{rho_al}
 2^{2s} \al^2 \leq \rho(s,\al)< 2 \cdot 2^{2s}\al^2.
\eeq 
Note that if the largest such $j$ is smaller than $2s +20$, the integral in (\ref{B_shifted}) is regarded as vanishing. Furthermore, there is no finite ``largest'' such $j$ precisely when $\al=0$; in this case we simply consider $\rho(s,\al)$ to be zero. Nevertheless, the statement (\ref{rho_al}) remains trivially true in this case (with the modification that the second inequality is not strict), so we retain this definition, as it is intuitively useful.

\subsection{Motivation for the double decomposition method}\label{sec_dd_motiv}
Let us denote by $\Bcal_{\lam,s}$ the operator defined for each $s \geq 0$ by the Fourier multiplier $B_\lam(s;\theta,\phi)$ given in (\ref{B_shifted}). Our goal is now to prove two types of estimates for the family $\Bcal_{\lam,s}$, namely $(\ell^2,\ell^2)$ estimates of the form (\ref{22}) on the line $\Re(\lam)=1$ and $(\ell^1, \ell^\infty)$ estimates of the form (\ref{1inf}) on the line $\Re(\lam)=-2/k$.
Here, as in the work leading to \cite{SW3}, certain difficulties will cause us to proceed with yet a third type of decomposition: a double decomposition of the major arcs.

Let us first identify the difficulties that arise, and then turn to mediating them. Understanding the multiplier $B_\lam(s;\theta,\phi)$ will ultimately come down to studying a multiplier of the form
\beq\label{L_al_be}
 L_\lam(\al,\be)  = \chi(c2^{2s}\al) \int_{\rho(s,\al)}^{2^{-2s-19}} \frac{y^{k\lam/2-1} e^{-2\pi Q^*(\be)/(y+i\al)}}{(y+i\al)^{\frac{k}{2}}} dy.
 \eeq
For $\Re(\lam)=1$, we would like to bound the multiplier $L_\lam$ directly, and for $\Re(\lam)=-2/k$, we would like to obtain a bound for the Fourier transform of $L_\lam$. To show, rather imprecisely, why it is difficult to obtain good bounds for $L_\lam$, note that taking the Fourier transform with respect to $\be$ alone (see equation (\ref{Q_Ftrans})) yields
\[ L_\lam \hat{^\be} (\al, \eta) = |A|^{1/2} \chi(c2^{2s}\al)\int_{2^{2s}\al^2}^{2^{-2s-19}} y^{k\lam/2-1} e^{-2\pi Q(\eta)(y+i\al)} dy.\]
Here we have used (\ref{rho_al}), so that $\rho(s,\al)\approx 2^{2s}\al^2$. To obtain a bounded function after taking the Fourier transform with respect to $\al$, it would be sufficient to see that $L_\lam \hat{^\be} (\al, \eta)$ is integrable with respect to $\al$. But for $\Re(\lam) = -2/k$, the above integral is approximately of size $\chi(c2^{2s}\al)\al^{-2}2^{-2s}$, and this is far from being integrable as a function of $\al$ near the origin. 

One way to remedy this situation is to replace the lower limit $\rho(s,\al)$ in (\ref{L_al_be}), which is quadratic in $\al$, by something linear in $\al$. One might then hope that the resulting Fourier transform in $\al$ would yield $\al^{-1 + i\gamma}$, where $\gamma = \Im(\lam) \neq 0$, in place of $\al^{-2}$, allowing one to proceed. 
Suppose, for example, we were able to reduce the situation to considering instead multipliers $L_\lam^r$ of the form 
\[ L^r_\lam(\al,\be)  = \chi(2^{r}\al)  \int_{2^{-r}\al}^{2^{-r+1}\al} \frac{y^{k\lam/2-1 } e^{-2\pi Q^*(\be)/(y+i\al)}}{(y+i\al)^{\frac{k}{2}}} dy,\]
where $r \geq 1$. Then for $\Re(\lam)=1$, we would have approximately 
\beq\label{line1} ||L_\lam^r||_{\ell^\infty} \leq \al^{-k/2}  \int_{2^{-r}\al}^{2^{-r+1}\al} y^{k\lam/2-1 } dy \approx \al^{-k/2}(2^{-r}\al)^{k/2} \approx 2^{-rk/2}.
\eeq
This gives decay, as desired. On the other boundary line, $\Re(\lam)=-2/k$, taking the Fourier transform in $\be$ first, we obtain 
\[ L^r_\lam \hat{^\be} (\al, \eta) =  \chi(2^r\al)  \int_{2^{-r}\al}^{2^{-r+1}\al} y^{k\lam/2-1} e^{-2\pi Q(\eta)(y+i\al)} dy.\]
We will show that one may disregard the contribution of the exponential factor with only minimal error, and hence this integral is well approximated, for $\lam = -2/k + i\gamma$, by
\[  \chi(2^r\al) \int_{2^{-r}\al}^{2^{-r+1}\al} y^{-2 + i\gamma k/2}dy \approx \chi(2^r\al) (2^{-r}\al)^{-1 + i\gam k /2}.\]
This is now integrable in $\al$ near the origin, but we have introduced a growth factor of size $2^r$. However, as we will see later, we could even accommodate growth of size $2^{2r}$ in this term and still obtain nontrivial decay via interpolation with the decay factor $2^{-rk/2}$ that will ultimately come from (\ref{line1}). We will accomplish the transition from multipliers of the form $L_\lam$ to multipliers of the form $L_\lam^r$ via a double decomposition of the major arcs,  which will allow us to ``slow down'' the approach of $\rho(s,\al)$ to zero.

\section{A double decomposition of the major arcs}\label{sec_double_decomp}

We define a double decomposition of the major arcs $\M_j(a,b;q)$ for each $j$ by setting
\beq\label{double_decomp}
\M_j (a,b;q) = \Union_{s=0}^{j/2-10} \Union_{r=0}^{j/2-s} \M_j^r(a,b;q).
\eeq
For $r \geq 1$ we define:
\[ \M_j^r(a,b;q) = \{(\theta, \phi): 2^{r-1}2^{-j} < |\theta - a/q| \leq 2^r 2^{-j}, |\phi_i-b_i/q| \leq (2q)^{-1}, i=1, \ldots, k \}.\]
For $r =0$ we define:
\[ \M_j^0(a,b;q) = \{(\theta, \phi): |\theta - a/q| \leq 2^{-j}, |\phi_i-b_i/q| \leq (2q)^{-1}, i=1, \ldots, k \}.\]
Thus we obtain a decomposition of the multiplier $M_\lam$ defined in (\ref{M_mult_dfn}), corresponding to the main term of $\Theta_Q$ supported on the major arcs, as
\beq\label{nu_double_decomp}
M_\lam (\theta, \phi) = \sum_{j=1}^\infty \sum_{s=0}^{j/2-10} \sum_{r=0}^{j/2-s} \nu^\lam_{j,r,s}(\theta, \phi).
\eeq
Precisely, for $1 \leq r \leq j/2-s$, we define
\begin{multline*} 
\nu^\lam_{j,r, s}(\theta, \phi) =   \sum_{2^s \leq q < 2^{s+1}} \sum_{\bstack{a=1}{(a,q)=1}}^q \sum_{b \in (\Z/q\Z)^k}
\frac{S_Q(a,b;q)}{q^k} \tchi(2^{-r}2^j|\theta -a/q|) \psi_q(\phi-b/q) \\
	 \cdot \; \int_{2^{-j}}^{2^{-j+1}} \frac{y^{k\lam/2-1} e^{-2\pi Q^*(\be) /(y+i\al)}}{(y+i\al)^{\frac{k}{2}}} dy.
\end{multline*}
Here $\tchi$ denotes the characteristic function of the interval $(1/2,1]$. For $r=0$, $\nu^\lam_{j,0,s}(\theta, \phi)$ is defined similarly, but with $\tchi$ replaced by $\chi_{[0,1]}$. 
The function $\psi_q(\be)$ is a smooth, even, non-negative cut-off function on $\R^k$ with $|\psi_q| \leq 1$, such that $\psi_q(\be)=1$ where $-1/2q \leq \be_i \leq 1/2q$, and vanishes outside the interval $-3/4q \leq \be_i \leq 3/4q$ in each coordinate; we furthermore choose $\psi_q$ such that for $\phi$ in the unit cube,
\beq\label{psi_q_sum}
 \sum_{b \in (\Z/q\Z)^k} \psi_q(\phi-b/q) =1.
 \eeq

We now wish to interchange the order of summation in $j$ and $s$. Effectively, since the sum over $0\leq s \leq j/2-10$ is vacuous unless $j \geq 20$, we have
\[ \sum_{j=1}^\infty \sum_{s=0}^{j/2-10} \sum_{r=0}^{j/2-s} = \sum_{s=0}^\infty \sum_{r=0}^{9} \sum_{j = 2s + 20}^\infty +  \sum_{s=0}^\infty \sum_{r=10}^{\infty} \sum_{j = 2r + 2s}^\infty. \]
Thus we may write
\beq\label{nu_sum}
 M_\lam(\theta, \phi) = \sum_{s=0}^\infty \sum_{r=0}^\infty \nu_{r,s}^\lam (\theta, \phi),
 \eeq
where
\beq\label{nu_rs}
 \nu_{r,s}^\lam = \sum_{j \geq \max(2s+20, 2s + 2r)} \nu_{j,r,s}^\lam.
 \eeq
Define for each $r,s, \lam$, the corresponding operator 
\[ (\Vcal_{r,s}^\lam f)\hat{\;} (\theta, \phi) = \nu_{r,s}^\lam (\theta, \phi) \hat{f}(\theta, \phi).\]
Our goal is to prove the following explicit version of the main proposition, Proposition \ref{key_J_prop}:
\begin{prop}\label{interp_bds_prop}
For each $r \geq 0, s\geq 0$, and for an appropriate analytic function $\A(\lambda)$, for $\Re(\lam)=1$,
\beq\label{A_lam1}
 ||\A(\lam) \Vcal_{r,s}^\lam f||_{\ell^2} \leq A 2^{-rk/2}2^{-sk/2} ||f||_{\ell^2},
 \eeq
and for $\Re(\lam) = -2/k$,
\beq\label{A_lam2}
 ||\A(\lam) \Vcal_{r,s}^\lam f||_{\ell^\infty} \leq A 2^{r}2^{2s} ||f||_{\ell^1}.
 \eeq
As a consequence, for all $2/(k+4) < \lam \leq 1$ and $1/p + 1/p'=1$ with $1/p' = 1/p - k(1-\lam)/(k+2)$,
\[ ||\A(\lam)\Vcal_{r,s}^\lam f||_{\ell^{p'}} \leq A2^{-\del_1(\lam)r -\del_2(\lam)s}||f||_{\ell^p}\] for some  $\del_1(\lam), \del_2(\lam) >0$. 
\end{prop}
This interpolation is sharp with respect to $s$, in the sense that given the two bounds (\ref{A_lam1}) and (\ref{A_lam2}), the interpolated exponents decay only in the range $2/(k+4) < \lam \leq 1$, and for no smaller $\lam$. 
 (However, for $r$ alone, as noted previously, it would be sufficient to prove a bound on $\Re(\lam) = -2/k$ with growth $2^{2r}$ in place of $2^r$.)
We will later obtain similar results for the multipliers incorporating the remainder term of the theta function on the major arcs, and the full theta function on the minor arcs.

\subsection{The main term supported on the major arcs: $r \geq 1$}\label{sec_maj_main}
We will prove Proposition \ref{interp_bds_prop} first for $r \geq 1$ and handle the case $r=0$ separately in Section \ref{sec_r_zero}. In fact, we will further reduce our consideration to the case $r \geq 10$; the cases $1 \leq r \leq 9$ are handled in exactly the same manner but with slightly different notation, due to the summation over $j \geq 2s +20$ for these cases. Thus fix $r \geq 10$. 
Consider
\begin{multline}\label{nu_rs'}
\nu^\lam_{r,s}(\theta, \phi) =   \sum_{2^s \leq q < 2^{s+1}} \sum_{\bstack{a=1}{(a,q)=1}}^q \sum_{b \in (\Z/q\Z)^k}
\sum_{j \geq 2s+2r} \frac{S_Q(a,b;q)}{q^k|A|^{1/2}} \tchi(2^{-r}2^j|\theta -a/q|) \\
	\cdot \psi_q(\phi-b/q)  \int_{2^{-j}}^{2^{-j+1}} \frac{y^{k\lam/2-1} e^{-2\pi Q^*(\be)/(y+i\al)}}{(y+i\al)^{\frac{k}{2}}} dy.
\end{multline}
For fixed $r,\theta$ and a fixed pair $a,q$, the cut-off factor $\tchi (2^{-r}2^j|\theta -a/q|)$ guarantees that only one $j$ in the sum will appear non-vacuously, namely that $j$ such that $2^{-j-1} < 2^{-r}|\theta - a/q| \leq 2^{-j}$.
 Furthermore, the supports of translates of $\tchi (2^{-r}2^j|\theta|)$ by distinct pairs $a_1/q_1$, $a_2/q_2$ are disjoint, for if not, we would have
\[  \frac{1}{4 \cdot 2^{2s}} < \frac{1}{q_1q_2} \leq \left|\frac{a_1}{q_1} - \frac{a_2}{q_2}\right| \leq \left| \theta - \frac{a_1}{q_1}\right| + \left| \theta - \frac{a_2}{q_2} \right|  \leq 2\cdot 2^r2^{-j},\]
which is impossible, since by assumption $j \geq 2s +2r$.

\subsubsection{Multiplier bounds on $\Re(\lam)=1$}
We will first bound $\nu^\lam_{r,s}$ on the line $\Re(\lam)=1$. 
Recall from Lemma \ref{Gauss_sum_lemmaQ} that for each triple $a,b,q$ arising in $\nu_{r,s}^\lam$, 
\beq\label{S_bdJ}
q^{-k}|S_Q(a,b;q)| \leq c_kq^{-k/2} \approx 2^{-sk/2}.
\eeq
Next we extract decay in $r$ from the integral in (\ref{nu_rs'}).
For fixed $\theta, a,q$ set $\al = \theta - a/q$ and let $j(\al)$ be the single $j$ that appears in the sum, i.e. such that
\beq\label{j_al_relation}
2^{-j(\al)-1} < 2^{-r}|\al| \leq  2^{-j(\al)}.
\eeq
Then the integral corresponding to $j(\al)$ is of magnitude
\beq\label{bound1}
 \left|   \int_{2^{-j(\al)}}^{2^{-j(\al)+1}} \frac{y^{k\lam/2-1} e^{-2\pi Q^*(\be)/(y+i\al)}}{(y+i\al)^{\frac{k}{2}}} dy \right| \leq |\al|^{-k/2} \int_{2^{-j(\al)}}^{2^{-j(\al)+1}} y^{k\Re(\lam)/2-1} dy.
 \eeq
As a result of (\ref{j_al_relation}), this  is bounded by
\[  |\al|^{-k/2} \int_{2^{-r}|\al|}^{2^{-r+2}|\al|} y^{k\Re(\lam)/2-1} dy= |\al|^{-k/2} O( (2^{-r}|\al|)^{k\Re(\lam)/2}) =O( 2^{-rk/2}),\]
for $\Re(\lam)=1$. Using the decay (\ref{S_bdJ}) from the Gauss sum (as well as the disjointness as $a,q$ vary of the supports of the $\tchi$ factors and  the condition (\ref{psi_q_sum}) on sums of $\psi_q$), we obtain for (\ref{nu_rs'}) the bound:
\beq\label{nu_rs_bd}
 |\nu_{r,s}^\lam(\theta,\phi)| \leq A2^{-rk/2}2^{-sk/2}.
 \eeq
 This proves (\ref{A_lam1}).

\subsubsection{Fourier coefficient bounds on $\Re(\lam)=-2/k$}
Next we consider $\nu_{r,s}^\lam$ on the line $\Re(\lam) = -2/k$. Here we would like to bound, for an appropriate analytic function $\A(\lam)$, the Fourier coefficients of $\A(\lam)\nu_{r,s}^\lam(\theta, \phi)$.  This requires a more intricate analysis.

The presence of the cut-off functions $\tchi$ and $\psi_q$ allows $\nu_{r,s}^\lam$ to be extended to a periodic function on all of $\R^{k+1}$, so that the Fourier coefficients are well-defined.
In general, let $c_l(g)$ represent the $l$-th Fourier coefficient of a function $g$ on $[0,1]^{k+1}$, where $l=(l_1,l_2) \in \Z \times \Z^k$.  First note that bounds for the Fourier coefficients are indeed sufficient to prove the bound (\ref{A_lam2}) for the operator $\Acal(\lam) \Vcal_{r,s}^\lam$. For indeed (ignoring $\Acal(\lam)$ for the moment),  
\begin{eqnarray*}
\Vcal_{r,s}^\lam f(n,t) &= &\int_{[0,1]^{k+1}} \nu_{r,s}^\lam(\theta,\phi) \hat{f}(\theta,\phi) e^{2\pi i \theta t} e^{2\pi i \phi \cdot n} d \theta d \phi \\
	& = & \int_{[0,1]^{k+1}} \nu_{r,s}^\lam(\theta,\phi) \sum_{u,v} f(u,v)e^{-2\pi i u\cdot \phi} e^{-2\pi i v \theta} e^{2\pi i \theta t} e^{2\pi i \phi \cdot n} d \theta d \phi \\
	& = &  \sum_{u,v} f(u,v) \int_{[0,1]^{k+1}} \nu_{r,s}^\lam(\theta,\phi)  e^{-2\pi i \theta (v-t)} e^{-2\pi i \phi \cdot (u-n)} d \theta d \phi \\
	& = & \sum_{u,v} f(u,v) c_{v-t,u-n}(\nu_{r,s}^\lam).
\end{eqnarray*}
Therefore, if $|c_{l_1,l_2}(\nu_{r,s}^\lam)| \leq B$ for all $(l_1,l_2)$, then $||\Vcal_{r,s}^\lam f||_{\ell^\infty} \leq B ||f||_{\ell^1}$. In particular, if we can show that all the Fourier coefficients of $\nu_{r,s}^\lam$ are $O(2^r 2^{2s})$, we will have proved (\ref{A_lam2}).

To compute $c_l(\nu_{r,s}^\lam)$:
\beq\label{Fc_orig}
 c_{l} (\nu_{r,s}^\lam ) =  \sum_{2^s \leq q < 2^{s+1}} \sum_{\bstack{a=1}{(a,q)=1}}^q \sum_{ b \in (\Z/q\Z)^k} \frac{ S_Q(a,b;q)}{q^k |A|^{1/2}} e^{-2\pi i l_1 a/q} e^{-2\pi i l_2 \cdot b/q} \cdot c_{l}(\psi_q(\phi) \mu_{r,s}^\lam(\theta,\phi)),
 \eeq
where 
\beq\label{mu_dfnrs}
 \mu_{r,s}^\lam(\theta,\phi) = \sum_{j=2r + 2s}^\infty \tchi(2^{-r}2^j|\theta|) \int_{2^{-j}}^{2^{-j+1}} \frac{y^{k\lam/2 -1} e^{-2\pi Q^*(\phi)/ (y+i\theta)}}{(y+i\theta)^{\frac{k}{2}}} dy.
\eeq
Note that by making the change of variables $\theta \mapsto \theta +a/q,$ $\phi \mapsto \phi+b/q$, we have reserved the arithmetic quantities $a,b,q$ from the Fourier coefficient of $\mu_{r,s}^\lam$, so as to exploit this in bounding the resulting exponential sums. Applying the bound (\ref{Gauss_avg_boundQ}) to the sum over $a,b$ of the Gauss sum and exponential factors, (\ref{Fc_orig}) becomes 
\beq\label{cl_bd}
 |c_l  (\nu_{r,s}^\lam ) | \leq \sum_{2^s \leq q < 2^{s+1}} q |c_l(\psi_q(\cdot )\mu_{r,s}^\lam(\cdot, \cdot))| \leq 2^{2s}\sup_{2^s \leq q <2^{s+1}}|c_l(\psi_q(\cdot )\mu_{r,s}^\lam(\cdot, \cdot))|.
 \eeq
Due to the small support of $\tchi$ and $\psi_q$, we may in fact compute the Fourier transform instead of Fourier coefficients. In particular,
\[ c_l(\psi_q(\cdot )\mu_{r,s}^\lam(\cdot, \cdot)) = \hat{\psi}_q(l_2) * \hat{\mu}_{r,s}^\lam(l_1,l_2),\]
and since $\psi_q$ is smooth, its Fourier transform is uniformly in $L^1$, so (\ref{cl_bd}) will imply (\ref{A_lam2}) if we can bound $ \hat{\mu}_{r,s}^\lam(l_1,l_2)$ by $O(2^r)$ uniformly in $l_1, l_2$.

Note that the following identity holds:
 \beq \label{Q_Ftrans}
 \int_{\R^k}e^{-2\pi Q^*(\phi) / (y+i\theta)}  e^{-2\pi i \phi \cdot \eta} d\phi    
 	=|A|^{1/2}(y+i\theta)^{k/2}e^{-2\pi Q(\eta)(y+i\theta)};
 \eeq
this is  obtained simply by diagonalizing the quadratic form $Q^*$ and computing the Fourier transform of the resulting product of Gaussians.
As a result, first taking the Fourier transform of $\mu_{r,s}^\lam$ with respect to $\phi$ alone gives
\beq\label{Ft_phi2} 
(\mu_{r,s}^\lam) \hat{^\phi}(\theta, \eta)  =  |A|^{1/2}\sum_{j=2r + 2s}^\infty \tchi(2^{-r}2^j|\theta|) \int_{2^{-j}}^{2^{-j+1}} y^{k\lam/2 -1} e^{-2\pi Q(\eta)(y+i\theta)} dy;
\eeq
notice the factor of $(y+i\theta)^{k/2}$ in the denominator of (\ref{mu_dfnrs}) has conveniently been cancelled. Interestingly, we have now effectively taken two Fourier transforms, once in the transformation law for $\Theta_Q$, and a second time here; thus we have passed from $Q$ to $Q^*$ and now back again to $Q$. 

Now we will show that we can disregard the exponential factor $e^{-2\pi Q(\eta)y}$ in (\ref{Ft_phi2}) with acceptable error, using the rule of thumb that an exponential factor $e^{-x}$ may be approximated by $x^{-1}$ if $x>1$ and by $1$ if $x \leq 1$. 
Thus in (\ref{Ft_phi2}), first suppose that $j$ is such that $2^{-j}Q(\eta) >1$. Then taking the Fourier transform with respect to $\theta$ of these terms in $(\mu_{r,s}^\lam)\hat{^\phi}(\theta,\eta)$, we obtain
\begin{multline}\label{Ft_theta}
 \left| \int_{-\infty}^\infty e^{-2\pi i \theta \xi} \sum_{\bstack{j}{2^{-j}Q(\eta) > 1}} \tchi(2^{-r}2^j|\theta|) \int_{2^{-j}}^{2^{-j+1}} y^{k\lam/2 -1} e^{-2\pi Q(\eta)(y+i\theta)} dy d\theta \right|  \\
 \leq  \sum_{\bstack{j}{2^{-j}Q(\eta)> 1}} \int_{2^{-j}2^{r-1} < |\theta| \leq  2^{-j}2^r} \int_{2^{-j}}^{2^{-j+1}} y^{k\Re(\lam)/2 -1} e^{-2\pi Q(\eta) y} dy d\theta. 
 \end{multline}
By the approximation $e^{-2\pi Q(\eta) y}  \leq (2\pi Q(\eta) y)^{-1}$, the inner integral is bounded by
\[ AQ(\eta)^{-1} \int_{2^{-j}}^{2^{-j+1}} y^{k\Re(\lam)/2 -2} dy \leq A2^{2j}Q(\eta)^{-1},\]
for $\Re(\lam) = -2/k$. Therefore (\ref{Ft_theta}) is bounded by
\beq\label{little_j}
A  \sum_{\bstack{j}{2^{-j}Q(\eta) > 1}}  2^{2j}Q(\eta)^{-1}\int_{2^{-j}2^{r-1} < |\theta| \leq  2^{-j}2^r}  d \theta \leq A 2^r Q(\eta)^{-1}\sum_{\bstack{j}{2^{-j}Q(\eta) > 1}}2^j \leq A 2^r.
  \eeq

We now consider the remaining $j$ in (\ref{Ft_phi2}), for which $2^{-j}Q(\eta)  \leq 1$. Then 
\[ |e^{-2\pi Q(\eta) y} - 1| \leq 2\pi Q(\eta) |y| ,\] 
so after taking the Fourier transform with respect to $\theta$, the error of replacing the exponential factor $e^{-2\pi Q(\eta) y}$ by $1$ in these terms in (\ref{Ft_phi2}) is of size
\beq\label{approx_error}
A \sum_{\bstack{j\geq 2s +2r}{2^{j}\geq Q(\eta)}} \int_{2^{-j}2^{r-1} < |\theta| \leq  2^{-j}2^r} \int_{2^{-j}}^{2^{-j+1}} y^{k\Re(\lam)/2}  Q(\eta) dy d\theta \leq A 2^r Q(\eta) \sum_{\bstack{j\geq 2s +2r}{2^{j} \geq Q(\eta) }} 2^{-j} \leq A2^r.    
 \eeq
Thus, finally, we are left with the terms in (\ref{Ft_phi2}) with $j$ such that $2^{-j} Q(\eta) \leq 1$, but with the exponential factor $e^{-2\pi Q(\eta) y}$ replaced by $1$. We compute the Fourier transform of these terms in $(\mu_{r,s}^\lam)\hat{^\phi}(\theta,\eta)$ with respect to $\theta$:
\begin{multline}\label{Ft_theta2}
\sum_{2^j \geq Q(\eta)} \int_{-\infty}^\infty e^{-2\pi i  \theta \xi} \tchi(2^{-r}2^j |\theta|)  \int_{2^{-j}}^{2^{-j+1}} y^{k\lam/2 -1} e^{-2\pi i Q(\eta) \theta} dy d\theta \\
 = 
A \sum_{2^j \geq Q(\eta)} 2^{-jk\lam/2} \int_{2^{-j} 2^{r-1} < |\theta| \leq 2^{-j} 2^r } e^{-2\pi i (\xi + Q(\eta)) \theta} d\theta \\
 =  
A2^r \sum_{2^j \geq Q(\eta)} 2^{-j} 2^{-jk\lam/2} \int_{1/2 < |\theta| \leq 1 } e^{-2\pi i 2^{-j}2^r  (\xi + Q(\eta)) \theta} d\theta .
\end{multline}
Now we consider two further cases: temporarily let $\xi '  = \xi + Q(\eta)$ and first suppose that $2^{-j}2^r|\xi'|>1$; then these $j$ contribute (by calculating the integral directly)
\[ A2^r \sum_{\bstack{2^j \geq Q(\eta)}{2^{-j}2^r|\xi'|>1}}  2^{-j} 2^{-jk\Re(\lam)/2} (2^{-j}2^r |\xi'|)^{-1} 
	\leq A |\xi'|^{-1}  \sum_{2^j <2^r|\xi'|} 2^{-jk\Re(\lam)/2} \leq  A2^r,\]
again using the fact that $\Re(\lam) = -2/k$.

For those $j$ such that $2^{-j}2^r|\xi'| \leq 1$, we replace $e^{-2\pi i 2^{-j}2^r  \xi' \theta}$ by $1$ in (\ref{Ft_theta2}), which introduces an error of size $O(2^{-j}2^r |\xi'| |\theta|)$, which contributes
\beq\label{little_j2}
 A2^r \sum_{2^j  \geq 2^r |\xi'|} 2^{-j}2^{-j k\Re(\lam)/2} 2^{-j}2^r |\xi'| \leq A2^{2r}|\xi'| \sum_{2^j \geq 2^r |\xi'|} 2^{-j} \leq A2^r. \eeq
Thus finally we are left with those $j$ in (\ref{Ft_theta2}) for which $2^{-j}2^r|\xi'| \leq 1,$ with the exponential replaced by $1$. These terms contribute
\beq\label{AA}
 A2^r \sum_{\bstack{2^j \geq Q(\eta)}{2^j \geq 2^r|\xi'|}} 2^{-j(1 + k\lam/2)}.
 \eeq
This sum behaves like $O(\frac{1}{2^{1+k\lam/2}-1})$, which has a pole as $\lam$ approaches $-2/k$. Thus we must introduce the analytic function $\A(\lam) = 2^{1 + k\lam/2} -1;$ multiplying (\ref{AA}) by $\A(\lam)$, it is then bounded by $A2^r $.

In conclusion, we have shown that (with the introduction of the analytic function $\A(\lam)$) all the terms in $\hat{\mu}^\lam_{r,s}$ are bounded by $A2^r$, and hence in (\ref{cl_bd}), $|c_l(\nu_{r,s}^\lam)| \leq A 2^r 2^{2s}$. 
Therefore, we may conclude that
\beq\label{concl_AV}
 ||\A(\lam) \Vcal_{r,s}^\lam f ||_{\ell^\infty} \leq A 2^r 2^{2s} ||f||_{\ell^1}
 \eeq
for $\Re(\lam) = -2/k$, which completes the proof of Proposition \ref{interp_bds_prop} for $r \geq 1$.

\subsection{The main term supported on the major arcs: $r =0$}\label{sec_r_zero}
The term $\nu_{0,s}^\lam$ must be handled separately due to subtle considerations near the origin. Recall that 
\begin{eqnarray}\label{nu_zero} 
\nu_{0,s}^\lam(\theta, \phi) & = & \sum_{j =2s+20}^\infty  \nu^\lam_{j,0, s}(\theta, \phi) \nonumber \\ 
& = & \sum_{2^s \leq q < 2^{s+1}} \sum_{\bstack{a=1}{(a,q)=1}}^q \sum_{b \in (\Z/q\Z)^k}
\frac{S_Q(a,b;q)}{q^k|A|^{1/2}} \sum_{j =2s+20}^\infty \chi_{[0,1]}(2^j|\theta -a/q|) \psi_q(\phi-b/q) \nonumber \\
	& &  \cdot \int_{2^{-j}}^{2^{-j+1}} \frac{y^{k\lam/2-1} e^{-2\pi Q^*(\beta)/(y+i\al)}}{(y+i\al)^{\frac{k}{2}}} dy \nonumber \\
	& = & \sum_{2^s \leq q < 2^{s+1}} \sum_{\bstack{a=1}{(a,q)=1}}^q \sum_{b \in (\Z/q\Z)^k}
\frac{S_Q(a,b;q)}{q^k|A|^{1/2}} \chi_{[0,1]}(2^{2s+19}|\theta -a/q|) \psi_q(\phi-b/q) \nonumber \\
	& &  \cdot \int_{r(\al)}^{2^{-2s-19}} \frac{y^{k\lam/2-1} e^{-2\pi Q^*(\beta)/(y+i\al)}}{(y+i\al)^{\frac{k}{2}}} dy.
\end{eqnarray}
Here $r(\al)$ plays the role that $\rho(s,\al)$ played previously,  i.e. $r(\al)$ denotes $2^{-J}$ where $J$ is the largest $j$ such that $|\al| = |\theta -a/q| \leq 2^{-j}$. Thus in particular we may say that $|\al| \leq r(\al) <2|\al|$. As before, we use the convention that $r(\al)=0$ if there is no finite largest such $j$ (i.e. precisely when $\al=0$, in which case this relation, with non-strict inequalities, holds vacuously). We also use the convention that the above integral vanishes if $r(\al) \geq 2^{-2s-19}$. 

We now prove the bounds of Proposition \ref{interp_bds_prop} for the operator $\Vcal_{0,s}^\lam$, defined by
\[ (\Vcal_{0,s}^\lam f)\hat{\;} (\theta,\phi) = \nu_{0,s}^\lam (\theta, \phi) \hat{f}(\theta,\phi).\]

\subsubsection{Multiplier bounds on $\Re(\lam)=1$}
We first consider the case $\Re(\lam) =1$. We will replace the factor $(y+i\al)^{-k/2}$ by $y^{-k/2}$ for $y \geq r(\al)$, using the expansion
\beq\label{expansion}
(y+i\al)^{-k/2} = y^{-k/2} + O( \frac{k\al }{2}y^{-k/2-1}).
\eeq
The first term leads to an integral in (\ref{nu_zero}) of the form
\beq\label{nu_zero_leading}
\int_{r(\al)}^{2^{-2s-19}} y^{k\lam/2 -k/2 -1} e^{-2\pi Q^*(\be)/(y+i\al)} dy.
\eeq
We now consider two cases:  for $Q^*(\be) >y$, we make the approximation
\[ |e^{-2\pi Q^*(\be)/(y+i\al)} | \leq e^{-2\pi Q^*(\be)/y} \leq (2\pi Q^*(\be)/y)^{-1}.\]
Then the contribution to the integral (\ref{nu_zero_leading}) from this portion is at most
\[ \int_0^{Q^*(\be)} y^{-1} ( y/Q^*(\be)) dy  = O(1).\]
If $y \geq Q^*(\be)$ then we approximate the exponential in (\ref{nu_zero_leading}) by replacing it by $1$. This introduces an error of size
\[ e^{-2\pi Q^*(\be)/(y+i\al)} -1 = e^{2\pi Q^*(\be) \al i/(y^2 + \al^2)}  \left( e^{-2\pi Q^*(\be) y/(y^2 + \al^2)}  - e^{-2\pi Q^*(\be)\al i/(y^2 + \al^2)} \right),\]
which in magnitude is
\beq\label{mag_exp}
O(Q^*(\be) ( |y/(y^2 + \al^2)| + |i\al/(y^2 + \al^2)|)) = O(Q^*(\be)(1/y+ \al/y^2)), 
\eeq
since trivially  $y^2 + \al^2 \geq y^2$. 
Now recall that in the range of the integral in (\ref{nu_zero_leading}), $y \geq r(\al) \geq |\al|$, thus (\ref{mag_exp}) is in fact $O(Q^*(\be)/y)$. Placing this estimate in the integral (\ref{nu_zero_leading}), we obtain a contribution from the error for $y \geq Q^*(\be)$  of at most 
\[ \int_{Q^*(\be)}^1 y^{-1} \cdot \frac{Q^*(\be)}{y} dy = O( 1 +  Q^*(\be)) = O(1).\]
Thus finally we consider the case where $y  \geq Q^*(\be) $ and we have replaced the exponential factor in (\ref{nu_zero_leading}) by $1$:
\[\left | \int_{Q^*(\be)}^1 y^{k\lam/2 - k/2 -1}dy \right|  = \left| \int_{Q^*(\be)}^1 y^{-1 + i\ga k/2} dy \right| \leq \frac{A}{|\gam|}, \]
where $\lam =1 + i\gamma$.
Putting all of these estimates together,  the integral in the multiplier $\nu_{0,s}^\lam$ resulting from the main term of the expansion (\ref{expansion}) for $(y+i\al)^{-k/2}$ is uniformly bounded by $ A(1/|\gam| + 1)$.

We now consider the contribution of the remainder term in the expansion (\ref{expansion}) (which is present only if $\al \neq 0$). This is of magnitude
\[  \int_{r(\al)}^{2^{-2s-19}} y^{k\Re(\lam)/2-1} O(\al  y^{-k/2-1}) dy  =  O(\al \int_{r(\al)}^{2^{-2s-19}} y^{-2}dy).\]
Recalling that $|\al| \leq r(\al) < 2|\al|$, this term is $O(\al ( 2^{2s+19} - \al^{-1}))$, and since the factor $ \chi_{[0,1]}(2^{2s+19}|\theta -a/q|) $ forces $|\al| 2^{2s+19}\leq 1$, this term is $O(1)$, uniformly in $\al$. 
(A more rigorous analysis may be performed for $\al$ arbitrarily close to the origin by redefining for every $\ep>0$ the arc closest to the origin by
\[ \M_{j,\ep}^0(a,b;q) = \{(\theta, \phi): \ep \leq |\theta - a/q| \leq 2^{-j}, |\phi_i-b_i/q| \leq (2q)^{-1}, i=1, \ldots, k \}.\] 
Repeating the above analysis for the multiplier $E_{\ep}$ corresponding to the remainder term $O(\al y^{-k/2-1}),$ supported on the arc $\M_{j,\ep}^0(a,b;q)$, again leads to an $O(1)$ bound, uniform in $\ep$.)

Thus in all cases we have bounded the integral  in (\ref{nu_zero})  by $O(1/|\ga| + 1)$.
We now use (\ref{Gauss_boundQ}) of Lemma \ref{Gauss_sum_lemmaQ} to bound the contribution of the Gauss sum, and by the disjointness of the major arcs we see that for $\Re(\lam)=1$,
\[ |\nu_{0,s}^\lam (\theta, \phi) | \leq A2^{-sk/2}(1/|\gam| +1).\]
Thus (choosing the analytic function $\Acal(\lam)$ to include the factor $1-\lam$ as well), the operator $\Vcal_{0,s}^\lam$ has the bound
\[ || \Acal(\lam)\Vcal_{0,s}^\lam f||_{\ell^2} \leq A_\ga2^{-sk/2}||f||_{\ell^2}\]
for $\Re(\lam)=1$.

\subsubsection{Fourier coefficient bounds on $\Re(\lam) =-2/k$}
We now consider the case $\Re(\lam)=-2/k$, for which we will bound the Fourier coefficients of $\A(\lam)\nu_{0,s}^\lam$ by $O(2^{2s})$, for an appropriate analytic function $\A(\lam)$. 
As before, because of the compact support of $\nu_{0,s}^\lam$ in $\theta$ and $\phi$, we may compute the Fourier transform directly; also as before we pull out the exponential factors $e^{-2\pi i l_1 a/q}$ and $e^{-2\pi i l_2 \cdot b/q}$ so as to include them in the Gauss sum. 
Recall from (\ref{Gauss_avg_boundQ}) of Lemma \ref{Gauss_sum_lemmaQ} that 
\beq\label{b_sumq}
\sum_{2^s \leq q < 2^{s+1}} \frac{1}{q^k} \left| \sum_{\bstack{a=1}{(a,q)=1}}^q \sum_{b \in (\Z/q\Z)^k}
S_Q(a,b;q)e^{-2\pi i l_1 a/q }e^{-2\pi i l_2 \cdot b/q} \right| \leq \sum_{2^s \leq q < 2^{s+1}} q = O(2^{2s}).
\eeq
This is the total growth we can allow in the Fourier coefficient $c_l(\nu_{0,s}^\lam)$, thus it remains to show that for every $l=(l_1,l_2)$, the remaining part of the Fourier coefficient, namely
$c_l(\psi_q(\cdot )\mu_{0,s}^\lam(\cdot, \cdot))$, is uniformly bounded,
where
\[ \mu_{0,s}^\lam(\theta, \phi) = \chi_{[0,1]}(2^{2s+19}|\theta|) \int_{r(\theta)}^{2^{-2s-19}}
\frac{y^{k\lam/2-1} e^{-Q^*(\phi)/(y+i\theta)}}{(y+i\theta)^{\frac{k}{2}}} dy.\]
Recall that $r(\theta)$ is $2^{-J}$ where $J$ is the largest $j$ such that $|\theta| \leq 2^{-j}$.

Since $\phi_q$ is a smooth function, it suffices to show that $\hat{\mu}_{0,s}^\lam$ is uniformly bounded.  Taking the Fourier transform first in $\phi$, we obtain (up to a constant)
\[ \mu_{0,s}^\lam \hat{^\phi} (\theta, \eta) = \chi_{[0,1]}(2^{2s+19}|\theta|) \int_{r(\theta)}^{2^{-2s-19}}
y^{k\lam/2-1} e^{-2\pi Q(\eta) (y+i\theta)} dy.\]
Then, taking the Fourier transform with respect to $\theta$, and using a dyadic decomposition of the resulting integral, 
\beq\label{mu_hat_0s}
 \hat{\mu}_{0,s}^\lam(\xi,\eta) =  \sum_{m = 2s+20}^\infty \int_{ \pm 2^{-m}}^{\pm 2^{-m+1}} e^{-2\pi i \theta (\xi +Q(\eta)  )} \int_{r(\theta)}^{2^{-2s-19}}
y^{k\lam/2-1} e^{-2\pi Q(\eta) y} dy d\theta.
\eeq
In particular, for each dyadic integral with $m$ fixed,  $2^{-m} \leq |\theta| \leq 2 \cdot 2^{-m}$ and hence $r(\theta) = 2 \cdot 2^{-m}$ in the innermost integral. This is very useful, as we can now switch orders of integration. 

We now consider two types of index $m$: first suppose that $2^m \leq |\xi'| $, where as before $\xi ' =\xi + Q(\eta)  $. These terms contribute
\beq\label{r0_m_sum}
  \sum_{\bstack{m \geq 2s+20}{2^m \leq |\xi'|}} \int_{2^{-m+1}}^{2^{-2s-19}}
y^{k\lam/2-1} e^{-2\pi Q(\eta) y} \int_{ \pm 2^{-m}}^{\pm 2^{-m+1}} e^{-2\pi i \theta \xi'}  d\theta dy.
\eeq
The innermost integral is bounded by $A/|\xi'|$ and thus when $\Re(\lam) = -2/k$, (\ref{r0_m_sum}) is of magnitude
\[ \frac{A}{|\xi'|}  \sum_{\bstack{m \geq 2s+20}{2^m \leq |\xi'|}} \int_{2^{-m+1}}^{2^{-2s-19}}
y^{-2} dy \leq \frac{A}{|\xi'|}  \sum_{2^m \leq |\xi'|} 2^m = O(1).\]
For $m$ such that $2^m > |\xi'| $, we will replace $e^{-2\pi i  \theta \xi'}$ by $1$, with an error of size $2\pi |\theta| |\xi'|$, which error gives a contribution to (\ref{mu_hat_0s}) of size
\[  A |\xi'|  \sum_{\bstack{m \geq 2s+20}{2^m > |\xi'|}} \int_{2^{-m+1}}^{2^{-2s-19}}
y^{-2} dy \int_{\pm 2^{-m}}^{\pm 2^{-m+1}} |\theta| d\theta \leq A|\xi'|  \sum_{2^m > |\xi'|} 2^{-m} = O(1).\]
Thus finally we consider $m$ such that  $2^m > |\xi'|$, but with $e^{-2\pi i \theta \xi'}$ replaced by $1$, which contribute to (\ref{mu_hat_0s}):
\beq\label{exp_one}
\sum_{\bstack{m \geq 2s+20}{2^m > |\xi'|}} 2\cdot 2^{-m}  \int_{2^{-m+1}}^{2^{-2s-19}}
y^{k\lam/2-1} e^{-2\pi Q(\eta) y} dy.
\eeq

We now further distinguish between two cases with respect to $Q(\eta) $ in (\ref{exp_one}). Suppose first that $2^{-m}Q(\eta) \geq1$. Then for $y$ in the range of the integral, we can make the approximation $|e^{-2\pi Q(\eta) y}| \leq (2\pi Q(\eta) y)^{-1}$, yielding
\[  A \sum_{\bstack{2^m \leq Q(\eta)}{2^m > |\xi'|}}  2^{-m}  \int_{2^{-m+1}}^{2^{-2s-19}}
y^{-2} \frac{1}{Q(\eta)y} dy \leq A Q(\eta)^{-1} \sum_{2^m \leq Q(\eta)} 2^m = O(1).\]

Next suppose that $2^{-m}Q(\eta) <1$. 
Then there exists $m_0 \leq m$  such that 
\[ 2^{-m_0} Q(\eta)  \leq 1< 2^{-m_0+1}Q(\eta),\]
and we break the integral in (\ref{exp_one}) into two parts, 
\beq\label{int_parts}
 \int_{2^{-m+1}}^{2^{-2s-19}}
y^{k\lam/2-1} e^{-2\pi Q(\eta) y} dy= \int_{2^{-m+1}}^{2^{-m_0}}+\int_{2^{-m_0}}^{2^{-2s-19}}.
\eeq
(If $2^{-2s-19} Q(\eta)<1$, then this decomposition is vacuous and we merely consider $\int_{2^{-m+1}}^{2^{-2s-19}}$ in place of the first integral.) In the first of the two integrals on the right hand side of (\ref{int_parts}) we approximate the exponential by $1$, making an error of size $O(Q(\eta)y)$, so that the total contribution of this  integral to (\ref{exp_one}) is 
\beq\label{first_type}
   \sum_{2^{-m} Q(\eta)<1}  2^{-m}  \int_{2^{-m+1}}^{2^{-m_0}}
y^{k\lam/2-1}dy + O( \sum_{2^{-m} Q(\eta)<1}  2^{-m}  \int_{2^{-m+1}}^{2^{-m_0}}
y^{-2} (Q(\eta)y) dy). \eeq
We reserve the evaluation of the first term in (\ref{first_type}) until later. Making the change of variables $u=Q(\eta)y$ in the error term of (\ref{first_type}), it is of size
\[Q(\eta) \sum_{2^{-m} Q(\eta)<1}  2^{-m}  \int_{Q(\eta)2^{-m+1}}^{Q(\eta)2^{-m_0}}
u^{-1} du .\]
But by assumption, $u \leq 1$ in the range of the integral, so that $u^{-1} \leq u^{-2+\del}$ for any $0< \del <1$.  Thus the error term in (\ref{first_type}) is bounded above by
\[ Q(\eta) \sum_{2^{-m}Q(\eta)<1} 2^{-m}(Q(\eta)2^{-m+1})^{-1+\del} 
\leq A Q(\eta)^\del \sum_{2^{-m}Q(\eta)<1} 2^{-m\del} = O(1).\]	

In the second integral on the right hand side of (\ref{int_parts}), $Q(\eta)y \geq 1$, so that we estimate the exponential by $(Q(\eta)y)^{-1}$ rather than by $1$; thus the contribution of this portion of (\ref{int_parts}) to (\ref{exp_one}) is of size
\[ \sum_{2^{-m}Q(\eta)<1} 2^{-m} \int_{2^{-m_0}}^{2^{-2s-19}} y^{-2}(Q(\eta)y)^{-1}dy 
\leq Q(\eta)^{-1}2^{2m_0} \sum_{2^{-m}Q(\eta)<1} 2^{-m} .\]
This is $O(2^{2m_0}Q(\eta)^{-2}) = O(1)$, since by construction $1 \leq 2^{m_0} Q(\eta)^{-1} <2$.

Finally, all that remains is the estimation of the first term in (\ref{first_type}), which contributes to (\ref{exp_one}) the quantity
\beq\label{inft_sum}
\sum_{2^m> Q(\eta)} 2^{-m}  \int_{2^{-m+1}}^{2^{-m_0}}
y^{k\lam/2-1} dy=A \sum_{2^m > Q(\eta)} (2^{-m})^{1+k\lam/2}.
\eeq
This last sum has a singularity at $\lam=-2/k$, and thus we take our analytic function to be $\A(\lam) = 2^{1 + k\lam/2} - 1.$ We may then conclude that $\A(\lam)\hat{\mu}_{0,s}^\lam$ is uniformly bounded as $\Re(\lam)$ approaches $-2/k$. 

Combining this with the result (\ref{b_sumq}) for the Gauss sum component of $c_l(\nu_{0,s}^\lam)$, we have bounded the Fourier coefficients of $\A(\lam) \nu_{0,s}^\lam$ by $O(2^{2s})$, and hence 
\[ ||\A(\lam) \Vcal_{0,s}^\lam f ||_{\ell^\infty} \leq A 2^{2s} ||f||_{\ell^1}\]
for $\Re(\lam) = -2/k$.
This completes the proof of Proposition \ref{interp_bds_prop} in the case $r=0$, and hence concludes our discussion of the multiplier corresponding to the main term of $\Theta_Q$ supported on the major arcs.

\section{The contribution of the remainder term on the major arcs}\label{sec_maj_rem}

It remains to consider two more portions of the multiplier $\nu_\lam(\theta,\phi)$ defined in (\ref{nu_dfn}):  that arising from the remainder term  $E_{\lam}(y+i\theta,\phi)$ of the theta function supported on the major arcs, and that arising from the theta function $\Theta_Q(y+i\theta,\phi)$ (with the main term and the remainder term taken together) supported on the minor arcs. For both of these remaining components we will not require the double decomposition of the major arcs introduced in Section \ref{sec_double_decomp}; instead, we merely use the original major/minor decomposition introduced in Section \ref{sec_mm_decomp}, with major arcs $\M_{j,s}(a,b;q)$ defined as in (\ref{Maj}). (We retain the smooth cut-off $\psi_q$ defined in (\ref{psi_q_sum}) to demarcate the arcs with respect to the variable $\phi$.)

We first consider the contribution of the remainder term $E_\lam(y+i\theta,\phi)$ in the approximate identity for the theta function $\Theta_Q$, which we recall is defined in (\ref{remainder_dfnQ}) by
\beq\label{remainder_dfn}
 E_\lam(y+i\theta,\phi) = \sum_{\bstack{m \in \Z^k}{m \neq 0}} \frac{S_Q(a,b-m;q)}{q^k|A|^{1/2}} \cdot \frac{e^{-2\pi Q^*(\be+m/q)/(y+i\al)}}{ (y+i\al)^{\frac{k}{2}}}.
 \eeq
We define $\Ecal_\lam$ to be the operator with Fourier multiplier corresponding to this term, supported on the major arcs. Recall that in the original major/minor decomposition of Section \ref{sec_mm_decomp}, for each fixed $j$ we performed a dissection of the unit interval of level $2^{j/2}$, and then defined the major arcs by taking $q \approx 2^s$ and letting $s$ range up to $s \leq j/2-10$. 
For each $j,s$ pair, let $\chi_{\M_{j,s}}(\theta,\phi)$ denote the characteristic function of the union of all the major arcs for that fixed $j,s$ pair.
Then let $\mathcal{E}_{\lam,j,s}$ be the operator with Fourier multiplier
\beq\label{remainder_mult}
 E_{\lam,j,s}(\theta, \phi) = \chi_{\M_{j,s}}(\theta,\phi) \int_{2^{-j}}^{2^{-j+1}} y^{k\lam/2-1}E_\lam(y+i\theta,\phi) dy.\eeq
This time we will sum first in $s$ and then in $j$, defining the remainder term operator by
\beq\label{remainder_sum}
\mathcal{E}_\lam f= \sum_{j=1}^\infty \mathcal{E}_{\lam,j}f = \sum_{j=1}^\infty \sum_{s \leq j/2-10}\mathcal{E}_{\lam,j,s} f ,
 \eeq
where $\mathcal{E}_{\lam,j}$ has multiplier $E_{\lam,j} = \sum_{s \leq j/2-10} E_{\lam,j,s}$.
We will prove the following analogue of the main Proposition \ref{key_J_prop}:
\begin{prop}\label{major_arcs_error_terms}
For every $j \geq 1$, for $\Re(\lam) = 1,$
\[ ||\mathcal{E}_{\lam,j}f||_{\ell^2} \leq A 2^{-jk/4} ||f||_{\ell^2},\]
and for $\Re(\lam) = -2/k,$
\[ ||\mathcal{E}_{\lam,j}f||_{\ell^\infty} \leq A 2^{j\al_k} ||f||_{\ell^1},\]
where $\al_1=1$, $\al_2 = 1 + \ep$ for any $\ep>0$, and $\al_k=1/2 + k/4$ for $k \geq 3$.
As a consequence, for all $\max(2/(k+4), k/(2k+2))< \lam \leq 1$ and $1/p + 1/p'=1$ with $1/p' = 1/p - k(1-\lam)/(k+2)$,
 \[ ||\mathcal{E}_{\lam,j}f||_{\ell^{p'}} \leq A 2^{-\del(\lam) j} ||f||_{\ell^p},\]
 for some $\del(\lam) >0$.
 \end{prop}

\subsection{Multiplier bounds on $\Re(\lam)=1$}
Recall from Proposition \ref{approx_id_no_derivQ} that for $2^{-j} \leq y \leq 2^{-j+1}$ and $\theta, \phi$ in a major arc associated to $a,b,q$,
\[
|E_\lam(y+i\theta,\phi)| =O(y^{-k/4}).\]
Applying this in $E_{\lam,j,s}$, 
\beq\label{E_bound'}
 |E_{\lam,j,s}(\theta,\phi) | \leq  A\chi_{\M_{j,s}} (\theta,\phi) \int_{2^{-j}}^{2^{-j+1}} y^{k\Re(\lam)/2-1} y^{-k/4} dy  \leq \chi_{\M_{j,s}}(\theta,\phi) 2^{- jk/4}.
 \eeq
Since $E_{\lam,j}$ is a finite sum over $s \leq j/2-10$, and since for a fixed $j$, the major arcs as $s$ varies are disjoint, we obtain the result of Proposition \ref{major_arcs_error_terms} for $\Re(\lam) = 1$.

\subsection{Fourier coefficient bounds on $\Re(\lam)=-2/k$}\label{sec_maj_remainder}
For $\Re(\lam)=-2/k$, it is sufficient to show that for all $l=(l_1,l_2)$ the Fourier coefficient
\[ |c_{l} ( E_{\lam,j}(\theta,\phi) )| \leq  A2^{j\al_k}.\]
We can compute the Fourier coefficient directly as 
\begin{multline}\label{Fc3}
 c_{l}(E_{\lam,j}(\theta,\phi))  = \sum_{s \leq j/2-10} \sum_{2^s \leq q < 2^{s+1}} \sum_{m \neq 0} \sum_{\bstack{a=1}{(a,q)=1}}^q \sum_{b \in (\Z/q\Z)^k} \frac{S_Q(a,b-m;q)}{q^k|A|^{1/2}} \\ \cdot e^{-2\pi i l_1 a/q} e^{-2\pi i l_2 \cdot b/q} c_{l}(W_{\lam,j,s,q,m}(\theta,\phi)),
 \end{multline}
where
\[ W_{\lam,j,s,q,m}(\theta,\phi) = \chi(2^{j/2}2^s \theta) \psi_q(\phi) 
\int_{2^{-j}}^{2^{-j+1}} \frac{y^{k\lam/2-1} e^{-2\pi Q^*(\phi + m/q)/(y+i\theta)}}{(y+i\theta)^{k/2}} dy.\]
We first consider the contribution from the Gauss sum. Summing first in $b$ shows that
\begin{eqnarray}
&& \sum_{b \modd{q}^k} \frac{S_Q(a,b-m;q)}{q^k} e^{-2\pi i l_1 a/q} e^{-2\pi i l_2 \cdot b/q} \label{avg_Gauss_b}\\
& = &  e^{-2\pi i l_1 a/q} \sum_{r \modd{q}^k} e^{-2\pi i Q(r) a/q} e^{2\pi i r \cdot m/q} ( q^{-k} \sum_{b \modd{q}^k} e^{-2\pi i r \cdot b/q} e^{-2\pi i l_2 \cdot b/q} ) \nonumber \\
& = & e^{-2\pi i l_1 a/q}  \sum_{r \modd{q}^k} e^{-2\pi i Q(r) a/q} e^{2\pi i r \cdot m/q} \del_q(r + l_2), \nonumber
\end{eqnarray}
where $\del_q(r+l_2) = 1$ if and only if $r + l_2 \con 0 \modd{q}^k$. Thus (\ref{avg_Gauss_b}) is $O(1)$ 
and therefore
\beq\label{W_sum}
 |c_{l}(E_{\lam,j}(\theta,\phi)) | \leq A \sum_{s \leq j/2-10} \sum_{2^s \leq q < 2^{s+1}} \sum_{\bstack{a=1}{(a,q)=1}}^q |c_{l}(\sum_{m \neq0}W_{\lam,j,s,q,m}(\theta,\phi))|.
 \eeq
 
We will bound the Fourier coefficient of $W$ trivially by the $L^1$ norm of $W$. We first note that
 \beq\label{boundW}
  \sum_{m \neq 0} |W_{\lam,j,s,q,m}| \leq A \chi(2^{j/2}2^s \theta)\psi_q(\phi) \int_{2^{-j}}^{2^{-j+1}} \frac{y^{k\Re(\lam)/2-1}}{|y+i\theta|^{k/2}} \left| \sum_{m\neq 0} e^{-2\pi Q^*(\phi +m/q)/(y+i\theta)}  \right| dy.
  \eeq
Bounding the sum over $m$ by $O(y^{-k/4}q^{k/2}(y^2 + \theta^2)^{k/4})$ as in (\ref{general_m_sum_bound}), 
\beq\label{SumW}
 \sum_{m \neq 0} |W_{\lam,j,s,q,m}(\theta,\phi)| \leq A \chi(2^{j/2} 2^s \theta)\psi_q(\phi) q^{k/2} 2^{-\frac{jk}{2}(\Re(\lam)-1/2)}.
 \eeq
Hence, using the small supports of the $\chi, \psi_q$ factors as well as the fact that $q \approx 2^s$,
\begin{eqnarray*}
 |c_{l}(\sum_{m \neq0}W_{\lam,j,s,q,m}(\theta,\phi))| 
 &\leq & || \sum_{m \neq 0} W_{\lam,j,s,q,m} ||_{L^1([0,1]^{k+1})} \\
 &\leq & A 2^{-j/2}2^{-s(k/2+1)} 2^{-\frac{jk}{2}(\Re(\lam)-1/2)}.
 \end{eqnarray*}
Finally, using this in (\ref{W_sum}) and summing trivially over $a,q$ and then over $s\leq j/2-10$, we obtain
\beq\label{Rbound}
 |c_{l}(E_{\lam,j}(\theta,\phi)) | \leq A2^{j\al_k},
 \eeq
 where $\al_1=1$, $\al_2 = 1 + \ep$ for any $\ep>0$, and $\al_k=1/2 + k/4$ for $k \geq 3$.
This completes the proof of Proposition \ref{major_arcs_error_terms}.

\section{The contribution of the minor arcs}\label{sec_min_full}
We now turn to the contribution to the multiplier $\nu_\lam(\theta,\phi)$ of the theta function supported on the minor arcs.
The minor arcs, the complement of the major arcs, are described by those $\theta$ such that if $|\theta -a/q| \leq 1/(q2^{j/2})$ for some $(a,q)=1$, then $2^{j/2-10} < q \leq 2^{j/2}.$ Recall from (\ref{nu_j_dfn}) that the original dyadic multiplier is defined
\[ \nu_{\lam,j}(\theta,\phi) =  \int_{2^{-j}}^{2^{-j+1}} y^{k\lam/2-1} \Theta_Q(y+i\theta,\phi) dy \]
where $\Theta_Q(y+i\theta,\phi)$ may be expressed as in (\ref{approx_id_eqn_JQ}) as the sum of a main term and a remainder term. Let $\Vcal^\lam_{\m,j}$ be the operator with Fourier multiplier $\chi_{\m_j}(\theta,\phi)\nu_{\lam,j}(\theta,\phi),$ where $\chi_{\m_j}$ represents the characteristic function of the union of the minor arcs for that fixed $j$. 
We will prove the final remaining case of the main Proposition \ref{key_J_prop}: 
\begin{prop}\label{V_prop}
 For every $j \geq 1$, for $\Re(\lam) =1$, 
\beq\label{minor_error_2}
 ||\Vcal^\lam_{\m,j} f||_{\ell^2} \leq A2^{-kj/4} ||f||_{\ell^2},
 \eeq
and for $\Re(\lam) = -2/k$,
\beq\label{minor_error_1}
 ||\Vcal^\lam_{\m,j} f||_{\ell^\infty} \leq A 2^{j\al_k} ||f||_{\ell^1},
 \eeq
 where $\al_1=1$, $\al_2 = 1 + \ep$ for any $\ep>0$, and $\al_k=1/2 + k/4$ for $k \geq 3$.
As a consequence, for all $\max(2/(k+4), k/(2k+2)) < \lam \leq 1$ and $1/p + 1/p'=1$ with $1/p' = 1/p - k(1-\lam)/(k+2)$, 
\[ ||\Vcal_{\m,j}^\lam f||_{\ell^{p'}} \leq A2^{-\del(\lam)j} ||f||_{\ell^p}\] 
for some $\del(\lam) >0$.
\end{prop}

\subsection{Multiplier bounds on $\Re(\lam)=1$}
For $\Re(\lam)=1$, we will bound $\Theta_Q$ by
\beq\label{F_bd_minor}
|\Theta_Q(y+i\theta,\phi)| \leq A \frac{2^{-kj/4}}{|y+i\al|^{k/2}},
\eeq 
where $y \approx 2^{-j}$,  $\theta = a/q + \al$, $\phi = b/q + \be$. In fact,
by the usual bound (\ref{Gauss_boundQ}) for the Gauss sum and the fact that in the minor arcs $q > 2^{j/2-10}$, it is immediate from (\ref{approx_id_eqn_JQ}) that the main term of $\Theta_Q$ is bounded by 
\[ Aq^{-k/2} |y+i\al|^{-k/2} \leq A 2^{-kj/4} |y+i\al|^{-k/2} .\]
Thus it suffices to note that by the usual bound (\ref{E_bdQ}), for $2^{-j} \leq y \leq 2^{-j+1}$, 
\[|E_\lam(y+i\theta,\phi)| =O(y^{-k/4}) = O(2^{jk/4}),\]
 which is smaller than (\ref{F_bd_minor}) since $y \approx 2^{-j}$ and $q \approx 2^{j/2}$.
From (\ref{F_bd_minor}) it follows immediately that for $\Re(\lam)=1$,
\beq\label{chi_bound}
 |\chi_{\m_j}(\theta,\phi) \nu_{\lam,j}(\theta,\phi)| \leq 
 A\chi_{\m_j}(\theta,\phi) 2^{-kj/4} \int_{2^{-j}}^{2^{-j+1}} \frac{y^{k/2-1}}{|y+i\al|^{k/2}}dy \leq A 2^{-kj/4} ,
 \eeq
which gives (\ref{minor_error_2}).

\subsection{Fourier coefficient bounds on $\Re(\lam)=-2/k$}
For the bound (\ref{minor_error_1}) we will show that $|c_{l}(\chi_{\m_j} \nu_{\lam,j})| \leq A2^{j\al_k}$ for all $l=(l_1,l_2)$. It is simplest to do this indirectly: let $\M_j$ denote the union of the major arcs for a fixed $j$, so that $\chi_{m_j} = 1- \chi_{\M_j}$, where $\chi_{\m_j}$ and $\chi_{\M_j}$ are the characteristic functions of the minor and major arcs for that fixed $j$, respectively. Then
\beq\label{indirect}
 |c_{l}(\chi_{\m_j} \nu_{\lam,j})| \leq |c_{l}( \nu_{\lam,j})| + |c_{l}(\chi_{\M_j} \nu_{\lam,j})|.
 \eeq
We will bound each of the terms on the right hand side of (\ref{indirect}) separately. It may seem counterintuitive that we no longer require a double decomposition to bound the major arc contribution $c_l(\chi_{\M_j} \nu_{\lam,j})$. However, note that in the present argument we accept a relatively large bound on the line $\Re(\lam)=-2/k$, because we have a good bound of size $2^{-kj/4}$ on the line $\Re(\lam)=1$. This good bound relies crucially on the fact that we are considering minor arcs, which allows us to use a lower bound on $q$ to obtain the bound (\ref{F_bd_minor}) for $\Theta_Q$; clearly such a lower bound does not hold in the major arcs.

We bound the first term on the right hand side of (\ref{indirect}) trivially.
Recall that by the definition of the theta function,
\[ \nu_{\lam,j}(\theta,\phi) =  \int_{2^{-j}}^{2^{-j+1}} y^{k\lam/2-1}\sum_{m \in \Z^k} e^{-2\pi Q(m)(y+i\theta)} e^{-2\pi i m \cdot \phi} dy, \]
in which the sum is uniformly and absolutely convergent for $y \geq 2^{-j}$.
Thus, explicitly
\beq\label{computeFT}
 c_{l} (\nu_{\lam,j}) =   \int_{2^{-j}}^{2^{-j+1}} \sum_{m \in \Z^k}  \int_{[0,1]^{k+1}} 
	 e^{-2\pi i \theta (Q(m) + l_1)} e^{-2\pi i \phi \cdot (m+l_2)}    d\theta d \phi \; y^{k\lam/2-1}e^{-2\pi Q(m)y} dy.
	 \eeq
Here the integrals over $\theta$ and $\phi$ yield $1$ if $Q(m)= - l_1$ and $m= - l_2$, and vanish otherwise. Thus only the term $m= - l_2$ remains in the infinite sum and so for $\Re(\lam)=-2/k$ and $-l_1 = Q(l_2)$,
\beq\label{result}
 |c_{l} ( \nu_{\lam,j} )| \leq \int_{2^{-j}}^{2^{-j+1}} y^{k\Re(\lam)/2-1}dy \leq A 2^j ;
 \eeq
the Fourier coefficient otherwise vanishes. 

We now turn to the contribution of the major arcs in (\ref{indirect}); here we use the approximate identity for $\Theta_Q$ as in Proposition \ref{approx_id_no_derivQ}. We first consider the contribution of the main term; this has an $l=(l_1,l_2)$ Fourier coefficient of
\beq\label{F_coeff}
\sum_{s \leq j/2-10} \sum_{2^s \leq q < 2^{s+1}} \sum_{\bstack{a=1}{(a,q)=1}}^q \sum_{b\in (\Z/q\Z)^k} \frac{S_Q(a,b;q)}{q^k|A|^{1/2}} e^{-2\pi i l_1 a/q} e^{-2\pi i l_2 \cdot b/q} c_{l} (\psi_q(\phi) \mu_\lam(\theta,\phi))
\eeq
where 
\[ \mu_\lam(\theta, \phi) = \chi(2^{j/2}2^s \theta) \int_{2^{-j}}^{2^{-j+1}} \frac{y^{k\lam/2 -1} e^{-2\pi Q^*(\phi) /(y+i\theta)}}{(y+i\theta)^{k/2}} dy.\]
Since $\chi, \psi_q$ restrict the support of the integrand to a unit cube, we may again take the Fourier transform rather than computing Fourier coefficients. Also, since $\psi_q$ is smooth, it suffices to bound $\hat{\mu}_\lam$. As usual, we proceed to take the Fourier transform in $\phi$ first, yielding 
\beq\label{Ftphi}
 \mu_\lam \hat{^\phi}(\theta, \eta) = 
 \chi(2^{j/2}2^s \theta) \int_{2^{-j}}^{2^{-j+1}} y^{k\lam/2 -1} e^{-2\pi Q(\eta)(y+i\theta)} dy.
 \eeq
Taking the Fourier transform now with respect to $\theta,$ 
\[ \hat{\mu}_\lam (\xi, \eta) = 
 \int_{-\infty}^\infty \chi(2^{j/2}2^s \theta) e^{-2\pi i \theta(Q(\eta) + \xi)} \int_{2^{-j}}^{2^{-j+1}} y^{k\lam/2 -1} e^{-2\pi  Q(\eta) y} dy,\] 
 and so on the line $\Re(\lam) = -2/k$,
 \beq\label{mu_hat}
  |\hat{\mu}_\lam (\xi,\eta)| \leq \int_{|\theta| \leq 2^{-j/2} 2^{-s}} d\theta \int_{2^{-j}}^{2^{-j+1}} y^{k\Re(\lam)/2 -1} dy = A2^{j/2} 2^{-s}.
  \eeq
 We now use this, along with the Gauss sum bound (\ref{Gauss_avg_boundQ}) for the sum over $a,b$ in (\ref{F_coeff}) to show that (\ref{F_coeff}) is bounded in absolute value by 
 \[ A\sum_{s \leq j/2-10} \sum_{2^s \leq q < 2^{s+1}} \frac{1}{q^k} \cdot q^{k+1} 2^{j/2}2^{-s} \leq A 2^{j/2}\sum_{s \leq j/2-10} 2^{-s} \sum_{2^s \leq q < 2^{s+1}} q \leq A2^{j},\]
for any $\ep>0$.

All that remains is to compute the contribution of the error term $E_{\lam}(y+i\theta,\phi)$ on the major arcs to the Fourier coefficient $c_l(\chi_{\mathfrak{M}_j} \nu_{\lam,j}).$ But this known to be $O(2^{j\al_k})$ by precisely the computation we already performed in Section \ref{sec_maj_remainder}.
 Thus we obtain the bound (\ref{minor_error_1}) for $|| \V_{\m,j} f||_{\ell^\infty}$ as desired. This completes our discussion of the minor arcs. 
Taken together, Propositions \ref{interp_bds_prop}, \ref{major_arcs_error_terms}, and \ref{V_prop} prove the main proposition, Proposition \ref{key_J_prop}, and this completes the proof of Theorem \ref{dis_JQ_thm}.

\section{Appendix: Necessary conditions}\label{nec_sec}
In this appendix we provide examples showing that conditions (i) and (ii) are necessary for the $(\ell^p, \ell^q)$ boundedness of the operators $J_{Q_1,Q_2,\lam}$ considered in Theorem \ref{dis_JQ_thm}; we again take advantage of the flexibility to choose $Q_2=Q_1$, which we now denote simply by $Q$. 

\subsection{Representations by quadratic forms}
We will require a simple result on the average number of representations of a positive integer by a quadratic form. Given a positive definite quadratic form $Q$ in $k$ variables with integer coefficients, 
let $r_{Q,k}(n)$ denote the number of representations of $n$ by $Q$, and set 
\[A_{Q,k}(N) = \sum_{n=1}^N r_{Q,k}(n).\]
\begin{lemma}
There exist positive constants $C_1, C_2$, dependent on $Q$, such that 
\beq\label{A_bdQ}
 C_1 N^{k/2} + O(N^{(k-1)/2}) \leq A_{Q,k}(N) \leq C_2N^{k/2} + O(N^{(k-1)/2}).
 \eeq
\end{lemma}
This is a simple consequence of the comparability of $Q$ to the generic form  $|\cdot|^2$, for which a more precise statement holds, namely:
\[
 A_{|\cdot|^2,k}(N) = \frac{\pi^{k/2}}{\Gamma(k/2+1)}N^{k/2} + O(N^{(k-1)/2}).
\]
This latter relation may be obtained simply by comparing the number of integer lattice points $m \in \Z^k$ with $|m| \leq N^{1/2}$ to the volume of the corresponding Euclidean ball of radius $N^{1/2}$ (see for example \cite{Wal}).

\subsection{Condition (ii)}
We now turn to the necessity of conditions (i) and (ii) in Theorem \ref{dis_JQ_thm}, beginning with the simpler condition (ii). Define $f(n,t)=1$ if $(n,t)=(0,0)$ and $f(n,t)=0$ otherwise. Then clearly $f \in \ell^p(\Z^{k+1})$ for all $1 \leq p \leq \infty$. On the other hand,
\[ ||J_\lam f(n,t)||^q_{q}  = \sum_{\bstack{n\in \Z^k, t \in \Z}{t=Q(n)}} Q(n)^{-k\lam q/2} = \lim_{T \maps \infty} \sum_{t \leq T} r_{Q,k}(t) t^{-k\lam q/2}.\]
 By partial summation and (\ref{A_bdQ}), this is $O(T^{\frac{k}{2}(1-\lam q)})$ and hence in the limit is finite if $1/q <\lam$. The condition $1/p > 1-\lam$ follows by taking adjoints. 

\subsection{Condition (i)}
Showing that condition (i) is necessary in Theorem \ref{dis_JQ_thm} is slightly more involved. 
Define $f(n,t) = |t|^{-\al} \chi(n/|t|^{1/2})$ for some $\al>0$ to be chosen later, where $\chi$ denotes the characteristic function of the ``square annulus'' $\{x \in \R^k: 1/2 < |x_j| <2, 1 \leq j \leq k\}$. Then
\beq\label{fpp}
 ||f||_{p}^p = \sum_{n,t} |f(n,t)|^p = \sum_t \frac{1}{|t|^{\al p}} \sum_{\bstack{n_1, \ldots, n_k}{\frac{1}{2}|t|^{1/2} \leq |n_j| \leq 2|t|^{1/2}}} \chi \left(\frac{n}{|t|^{1/2}} \right)^p \approx c \sum_t \frac{|t|^{k/2}}{|t|^{\al p}}, 
 \eeq
which is finite if $\al >\frac{k+2}{2p}$.
Now consider 
\beq\label{J23}
 J_\lam f(n,t) = \sum_{\bstack{m \in \Z^k}{m \neq 0}} \frac{f(n-m, t-Q(m))}{Q(m)^{k\lam /2}}  \geq \sum_{\bstack{m_1, \ldots, m_k}{Q(m) \leq \del |t|}} \frac{1}{|t-Q(m)|^\al} \chi \left( \frac{n-m}{|t-Q(m)|^{1/2}} \right) \frac{1}{Q(m)^{k\lam /2}};
 \eeq
including the additional restriction $Q(m) \leq \del |t|$ for some small $\del >0$ gives a lower bound because all the summands are non-negative. As a result of this restriction, $t-Q(m) \approx t$ and so the last sum in (\ref{J23}) is approximately of size 
\[  \frac{1}{|t|^{\al}} \chi(\frac{n}{|t|^{1/2}}) \sum_{Q(m) \leq \del |t|} \frac{1}{Q(m)^{k\lam/2}} =  \frac{1}{|t|^{\al}} \chi \left(\frac{n}{|t|^{1/2}} \right) \sum_{l=1}^{\del |t|} \frac{r_{Q,k}(l)}{l^{k\lam /2}}  .\]
Applying partial summation and (\ref{A_bdQ}) shows that up to a constant,
\beq\label{J_lower_bd}
J_\lam f(n,t) \geq \frac{1}{|t|^\al} \chi \left(\frac{n}{|t|^{1/2}} \right) |t|^{\frac{k}{2}(1-\lam)}.
\eeq
Therefore, 
\[ ||J_\lam f||^q_{q} \geq \sum_t |t|^{\frac{kq}{2}(1-\lam)-\al q} \sum_n \chi \left(\frac{n}{|t|^{1/2}} \right)^q \approx \sum_t |t|^{\frac{kq}{2}(1-\lam)-\al q} |t|^{k/2}.\]
This last sum is finite if and only if $\al - \frac{k+2}{2q} > \frac{k}{2}(1-\lam).$ Recall that for $f$ to be in $\ell^p$ we required $\al>\frac{k+2}{2p}$; thus set $\al = \frac{k+2}{2p} + \ep$ for any $\ep >0$. Then to have $J_\lam f \in \ell^q$ as well, we must have $ \frac{k+2}{2p} - \frac{k+2}{2q} > \frac{k}{2}(1-\lam) -\ep$, and hence $1/p - 1/q \geq \frac{k(1-\lam)}{k+2}.$ This proves the necessity of condition (i) in Theorem \ref{dis_JQ_thm}.

\section*{Acknowledgments}
The author would like to thank Elias M. Stein for suggesting this area of inquiry, and for his generous advice and encouragement during many discussions of his previous work (with Stephen Wainger) on related problems. The author is also grateful to Michael Christ for providing a preprint of \cite{Chr88}.

\bibliographystyle{amsplain}
\bibliography{AnalysisBibliography}

\end{document}